\documentclass[a4paper,review]{cas-sc}

\usepackage[section]{placeins} 
\usepackage{amssymb}
\usepackage{amsthm}
\usepackage{amsmath}
\usepackage{upgreek}
\usepackage{pdflscape}
\usepackage{listings}
\usepackage{multirow}
\usepackage[percent]{overpic}
\usepackage{multicol}
\usepackage{color}
\usepackage{stmaryrd}
\usepackage{xfrac}
\usepackage[capitalise]{cleveref}
\usepackage{booktabs}
\usepackage[section]{placeins} 
\usepackage[section]{algorithm}
\usepackage{calc}
\usepackage[titletoc]{appendix}
\usepackage{siunitx}
\usepackage{algorithmic}
\usepackage[sort&compress,square,numbers]{natbib}

\usepackage{graphicx}
\usepackage{graphics}
\usepackage{tabularx}
\usepackage{wrapfig}
\usepackage{float}
\usepackage{subfig}
\usepackage[percent]{overpic}
\usepackage{varwidth}
\usepackage{tikz}
\usepackage{pgfplots}
\usepackage{adjustbox}
\usetikzlibrary{arrows,matrix,positioning,fit}
\usetikzlibrary{shapes,positioning}
\usetikzlibrary{backgrounds}
\usepackage{tikz-layers}
\usepgfplotslibrary{fillbetween}
\usetikzlibrary{intersections}
\pgfplotsset{compat=1.14}
\usepgfplotslibrary{colorbrewer}
\usepgfplotslibrary{patchplots}
\usepgfplotslibrary[colorbrewer]
\usetikzlibrary{pgfplots.colorbrewer}
\usetikzlibrary[pgfplots.colorbrewer]
\usepgfplotslibrary{units}
\usetikzlibrary{spy}
\usepackage{pgfplotstable}
\usepackage{arrayjobx}
\graphicspath{ {./figs/} }
\usetikzlibrary{external}
\usetikzlibrary{patterns}

\newlength\myheight
\newlength\mydepth
\settototalheight\myheight{Xygp}
\newcommand*\inlinegraphics[1]{%
  \settototalheight\myheight{Xygp}%
  \settodepth\mydepth{Xygp}%
  \raisebox{-\mydepth}{\includegraphics[height=\myheight]{#1}}%
}
\newcommand\orcid[1]{\href{https://orcid.org/#1}{\inlinegraphics{orcid_16x16.png}}}

\makeatletter
\def\BState{\State\hskip-\ALG@thistlm}
\makeatother

\newdefinition{definition}{Definition}[section]

\begin{document}

\title[mode=title]{Positivity-preserving entropy filtering for the ideal magnetohydrodynamics equations}
\shorttitle{Positivity-preserving entropy filtering for the ideal MHD equations}
\shortauthors{T. Dzanic \textit{et al.}}

\author[1]{T. Dzanic}[orcid=0000-0003-3791-1134]
\cormark[1]
\cortext[cor1]{Corresponding author}
\ead{tdzanic@tamu.edu}
\author[1]{F. D. Witherden}[orcid=0000-0003-2343-412X]

\address[1]{Department of Ocean Engineering, Texas A\&M University, College Station, TX 77843}

\begin{abstract}
In this work, we present a positivity-preserving adaptive filtering approach for discontinuous spectral element approximations of the ideal magnetohydrodynamics equations. This approach combines the entropy filtering method (Dzanic and Witherden, \textit{J. Comput. Phys.}, 468, 2022) for shock capturing in gas dynamics along with the eight-wave method for enforcing a divergence-free magnetic field. Due to the inclusion of non-conservative source terms, an operator-splitting approach is introduced to ensure that the positivity and entropy constraints remain satisfied by the discrete solution. Furthermore, a computationally efficient algorithm for solving the optimization process for this nonlinear filtering approach is presented. The resulting scheme can robustly resolve strong discontinuities on general unstructured grids without tunable parameters while recovering high-order accuracy for smooth solutions. The efficacy of the scheme is shown in numerical experiments on various problems including extremely magnetized blast waves and three-dimensional magnetohydrodynamic instabilities. 
\end{abstract}

\begin{keywords}
Discontinuous spectral element \sep
Ideal magnetohydrodynamics \sep
Shock capturing \sep
Positivity-preserving \sep
Entropy filtering
\end{keywords}



\maketitle

\section{Introduction}
\label{sec:intro}
The transport and interaction of a non-resistive conducting fluid and its electromagnetic field remain extensively investigated phenomena as they are instrumental in various applications ranging from the study of astrophysical accretion disks \citep{Hawley2000} and supernova remnants \citep{Hill2012} to magnetic confinement fusion \citep{Fasoli2016} and plasma physics \citep{Zhai2014}. These strongly nonlinear effects are governed by the equations of ideal magnetohydrodynamics (MHD), which are composed of a combination of the Euler equations of gas dynamics and Maxwell's equations of electromagnetism. From this formulation, a strong coupling between the magnetic field and the conducting fluid can be observed, where the magnetic field induces a current in the fluid which, in turn, gives rise to a second, induced magnetic field. This interaction can introduce multi-scale, multi-physics behavior in the system, such that magnetohydrodynamic flows can become exceedingly complex.

As a result of this complexity, the robust and accurate numerical approximation of ideal MHD can present many challenges. Since hyperbolic systems are known to produce discontinuities even with smooth initial conditions \citep{Hopf1950}, the numerical scheme must be able to robustly resolve these discontinuities which, in the case of MHD, come in the form of hydrodynamic and magnetic shocks and contact waves. Furthermore, the approximation of the ideal MHD equations also requires an intrinsic constraint on the solution in the form of a solenoidal magnetic field which may not be satisfied by the scheme even if the magnetic field is initially solenoidal. Without a mechanism to enforce this constraint, unphysical dynamics can arise in the solution which can lead to numerical instabilities. The standard numerical schemes for approximating MHD flows are finite difference and finite volume methods, whose properties and robustness are well-established in the literature \citep{Evans1988, Powell1999, Brio1988, Dai1998, Balsara1999, Balbs2006}. However, they possess certain drawbacks in that they are either difficult to extend to complex domains with unstructured grids or cannot recover high-order accuracy in a computationally efficient manner.

A particular class of schemes which have more recently grown in popularity are high-order discontinuous spectral element methods (DSEM) as they possess the geometric flexibility of finite volume methods while retaining the arbitrarily high-order accuracy and efficiency of spectral methods. As such, they provide a promising avenue for significantly decreasing the computational cost and expanding the viability of simulating complex MHD problems. However, due to the presence of discontinuities in MHD, DSEM approximations of these systems may introduce spurious oscillations in the solution in the form of Gibbs phenomena. Without proper treatment, these oscillations can result in unphysical predictions or the failure of the numerical scheme altogether. To extend to use of DSEM to MHD, various numerical stabilization techniques have been proposed, ranging from artificial viscosity methods \citep{Ciuc2020,Dzanic2022b} to limiting-type approaches \citep{RuedaRamrez2022, Wu2018}. While these various methods may be sufficient to stabilize the solution in many cases, they may not guarantee that the solution will abide by physical constraints, may require problem-dependent tunable parameters, can be computationally inefficient for general unstructured grids, or may be excessively dissipative in smooth regions of the flow. 

There is significant interest in the design of numerical schemes that are ``provably robust'' in the sense that they can guarantee that the solution will abide by certain physical constraints even in the presence of features such as discontinuities, the quintessential examples being positivity-preserving schemes for gas dynamics which guarantee the positivity of the density and internal energy/pressure. For DSEM, this property is typically achieved through some form of nonlinear limiting or filtering \citep{Zhang2011, Dzanic2022, Lin2023, RuedaRamrez2022}. However, designing schemes that possess this property without sacrificing the computational efficiency of DSEM for general unstructured grids and their advantageous scale-resolving properties in smooth flow regions can be challenging. In the context of MHD, this becomes even more difficult due to the additional complexity of the governing equations as well as the incorporation of differential constraints, namely solenoidal magnetic fields. As such, there is a need for numerical stabilization techniques for DSEM approximations of the ideal MHD equations that retain as many of these desirable properties as possible, namely that they: 
\begin{enumerate}
    \item Guarantee that physical constraints of the solution are satisfied.
    \item Are compatible with numerical techniques for enforcing intrinsic constraints such as a solenoidal magnetic field.
    \item Do not require problem-dependent tunable parameters.
    \item Do not appreciably degrade the ability of the underlying DSEM to resolve smooth portions of the flow.
    \item Can be easily and efficiently implemented on general unstructured grids.
\end{enumerate}

In this work, we propose a nonlinear adaptive filtering approach as a numerical stabilization technique for DSEM approximations of the ideal MHD equations to address these points. The proposed technique can be considered as an extension of the entropy filtering approach originally introduced by the authors for shock capturing in gas dynamics to the ideal MHD system \citep{Dzanic2022}. This technique relies on using the solution's ability to preserve convex invariants of the system, namely positivity of the density and pressure and a discrete local minimum entropy principle, to compute the necessary filter strength to ensure a well-behaved solution in the vicinity of discontinuities. Extending this approach to the ideal MHD system presents several challenges, primarily stemming from the treatment of the divergence-free constraint on the magnetic field. We utilize the eight-wave method of \citet{Powell1999} which introduces non-conservative source terms in the equation proportional to the divergence of the magnetic field. As these non-conservative terms can conflict with the necessary assumptions of the entropy filtering approach, we present a modified set of conditions and introduce an operator splitting approach to the system which allows the filtering method to retain its positivity-preserving properties for most practical applications. Furthermore, as the original approach for performing the optimization process necessary in the adaptive filtering framework as presented in \citet{Dzanic2022} was found to be quite computationally expensive, we develop a highly-efficient numerical approach which drastically reduces the overall computational cost. The resulting approach can robustly resolve strong hydrodynamic and magnetic discontinuities in the flow without appreciably degrading the accuracy of the underlying DSEM for smooth flows, does not require problem-dependent tunable parameters, and can be easily extended to unstructured grids with relatively low computational cost. The efficacy of the proposed method is demonstrated in a variety of numerical experiments including smooth transport, extremely magnetized blast waves, and three-dimensional magnetohydrodynamic instabilities computing using high-order approximations on both structured and unstructured grids.

The organization of this work is as follows. We present some preliminaries regarding DSEM approximations and the ideal MHD equations in \cref{sec:preliminaries}. The entropy filtering approach for ideal MHD is then introduced in \cref{sec:methodology}, and its numerical implementation and computational optimizations are presented in \cref{sec:implementation}. Results for various test cases are then shown in \cref{sec:results}, and conclusions are drawn in \cref{sec:conclusion}.
\section{Preliminaries}
\label{sec:preliminaries}
\subsection{Ideal magnetohydrodynamics}
The governing equations for the evolution of an ideal magnetohydrodynamic fluid can be given in the form of a hyperbolic conservation law as
\begin{equation}\label{eq:hyp_source}
        \partial_t \mathbf{u} + \boldsymbol{\nabla}{\cdot} \mathbf{F}\left(\mathbf{u}\right) = \mathbf{S}_{\boldsymbol{B}} (\mathbf{u}),
    \end{equation}
where $\mathbf{u} = \mathbf{u} (\mathbf{x}, t) \in \mathbb R^m$ is the solution of some number of field variables $m$ defined over a $d$-dimensional spatial domain $\mathbf{x} \in \mathbb R^d$ and time $t$, $\mathbf{F}(\mathbf{u}) \in \mathbb R^{m \times d}$, and $\mathbf{S_B}(\mathbf{u})$ is an additional source term to be defined in \cref{ssec:powell} whose purpose is to ensure a solenoidal magnetic field. The solution and flux are given as
    \begin{equation}
        \mathbf{u} = \begin{bmatrix}
                \rho \\ \boldsymbol{\rho}\mathbf{v} \\ \boldsymbol{B} \\ E
            \end{bmatrix} \quad  \mathrm{and} \quad \mathbf{F} = \begin{bmatrix}
                \boldsymbol{\rho}\mathbf{v}\\
                \boldsymbol{\rho}\mathbf{v}\otimes\mathbf{v} + \mathbf{I}\left(P + \frac{1}{2} \mathbf{B}{\cdot} \mathbf{B}\right) - \mathbf{B}\otimes\mathbf{B}\\
                \boldsymbol{v}\otimes\mathbf{B} - \mathbf{B}\otimes\boldsymbol{v}\\
            \left(E + P + \frac{1}{2} \mathbf{B}{\cdot} \mathbf{B}\right)\mathbf{v} - \mathbf{B}(\mathbf{v}{\cdot} \mathbf{B})
        \end{bmatrix},
    \end{equation}
where $\rho$ is the density, $\boldsymbol{\rho}\mathbf{v}$ is the momentum, $E$ is the total energy, $P = (\gamma-1)\left(E - \frac{1}{2}\rho\mathbf{v}{\cdot}\mathbf{v} - \frac{1}{2}\mathbf{B}{\cdot}\mathbf{B}\right)$ is the pressure, $\mathbf{B}$ is the magnetic field, and $\gamma$ is the specific heat ratio. Furthermore, the symbol $\mathbf{I}$ denotes the identity matrix in $\mathbb{R}^{d\times d}$ and $\mathbf{v} = \boldsymbol{\rho}\mathbf{v}/\rho$ denotes the velocity. The solution can be more conveniently expressed in terms of a vector of primitive variables as $\mathbf{q}=[\rho,\mathbf{v},\mathbf{B},P]^T$, and auxiliary quantities representing the magnetic pressure and plasma-beta can be defined as $P_b = \frac{1}{2}(\gamma - 1)\mathbf{B}{\cdot}\mathbf{B}$ and $\beta = 2P/(\mathbf{B}{\cdot}\mathbf{B})$, respectively. 

Due to the lack of magnetic monopoles, the MHD equations have an intrinsic constraint on the solution in the form of a solenoidal magnetic field, i.e.,
\begin{equation}
    \boldsymbol{\nabla}{\cdot}\mathbf{B} = 0.
\end{equation}
Although this constraint must be satisfied analytically by the MHD equations, numerical approximations do not necessarily satisfy it even if the magnetic field is initially solenoidal. If this constraint is not enforced by the scheme, numerical instabilities may arise in addition to the non-physical nature of the approximation. Many approaches exist to enforce this condition on the magnetic field, including the use of solenoidal basis functions \citep{Li2005}, projection methods \citep{Brackbill1980}, constrained-transport schemes \citep{Evans1988}, divergence cleaning methods \citep{Dedner2002}, and the eight-wave method \citep{Powell1999}. An overview of the salient techniques is presented in \citet{Wu2018}.

The entropy solution of \cref{eq:hyp_source} satisfies an entropy inequality of the form
\begin{equation}
    \partial_t\sigma(\mathbf{u}) + \boldsymbol{\nabla}{\cdot}\boldsymbol{\Sigma}(\mathbf{u}) \geq 0,
\end{equation}
where $(\sigma, \boldsymbol{\Sigma})$ is any numerical entropy-flux pair \citep{Harten1987} that satisfies the relation
\begin{equation*}
    \partial_{\mathbf{u}}\boldsymbol{\Sigma} = \partial_{\mathbf{u}} \sigma \partial_{\mathbf{u}} \mathbf{F}.
\end{equation*}
Note that this inequality may be negated depending on which notation is used for the numerical entropy. In \citet{Dao2022}, it was shown that the entropy solution (in a vanishing viscosity sense) of the ideal MHD system satisfies a minimum entropy principle on the specific physical entropy $\sigma = P\rho^{-\gamma}$ in the form
\begin{equation}
    \sigma \left (\mathbf{u}(\mathbf{x}, t + \Delta t )\right) \geq \underset{\mathbf{x}}{\min}\ \sigma \left (\mathbf{u}(\mathbf{x}, t)\right),
\end{equation}
where $\Delta t > 0$. This property is identical to the minimum entropy principle in gas dynamics \citep{Tadmor1986}, and it should be satisfied by the solution in both smooth regions and in the vicinity of discontinuities.

\subsection{Discontinuous spectral element methods}
 For nodal discontinuous spectral element approximations of \cref{eq:hyp_source}, including discontinuous Galerkin \citep{Hesthaven2008DG} and flux reconstruction \citep{Huynh2007} schemes, the domain $\Omega$ is partitioned into $N_e$ elements $\Omega_k$ such that $\Omega = \bigcup_{N_e}\Omega_k$ and $\Omega_i\cap\Omega_j=\emptyset$ for $i\neq j$. With a slight abuse of notation, the solution $\mathbf{u} (\mathbf{x})$ within each element $\Omega_k$ is approximated in a nodal manner as 
\begin{equation}
    \mathbf{u} (\mathbf{x}) = \sum_{i\in S} \mathbf{u}_i \phi_i(\mathbf{x}),
\end{equation}
where $\mathbf{x}_i\ \forall \ i \in S$ is a set of solution nodes, $\phi_i(\mathbf{x})$ are their associated nodal basis functions that possess the property $\phi_i(\mathbf{x}_j) = \delta_{ij}$, and $S$ is the set of nodal indices for the stencil. For brevity, we utilize the notation that $\mathbf{u}_i = \mathbf{u}(\mathbf{x}_i)$. The order of the approximation of the solution is denoted as $\mathbb P_p$ for some order $p$, where $p$ is the maximal order of $\mathbf{u} (\mathbf{x})$. This approximation formally yields a convergence rate of at least $p+1$ \citep{Hesthaven2008DG}.

The flux is approximated via the contribution of an interior term, denoted by the subscript $\Omega_k$, and an interface term, denoted by the subscript $\partial \Omega_k$, as 
\begin{equation}
    {\mathbf{F}}(\mathbf{u}) \approx {\mathbf{F}}_{\Omega_k}(\mathbf{u}) + {\mathbf{F}}_{\partial \Omega_k}(\mathbf{u}).
\end{equation}
For the interior component, the flux is computed through a collocation approach as 
\begin{equation}
    \mathbf{F}_{\Omega_k} (\mathbf{u} ) = \sum_{i\in S} \mathbf{F}(\mathbf{u}_i) \phi_i(\mathbf{x}),
\end{equation}
such that the interior contribution to the divergence of the flux can be computed as
\begin{equation}\label{eq:cij}
    \boldsymbol{\nabla}{\cdot}{\mathbf{F}}_{\Omega_k}(\mathbf{u}_i) = \sum_{j \in S} \mathbf{c}_{ij} \mathbf{F}(\mathbf{u}_j), \quad \mathrm{where} \quad \mathbf{c}_{ij} = \nabla \phi_i(\mathbf{x}_j).
\end{equation}
The interface component of the flux is formed over a set of interface nodes $\mathbf{x}_i \in \partial \Omega_k \ \forall \ i \in I$, where $I$ is a set of nodal indices for the interface stencil. We assume that these interface nodes are a subset of the solution nodes (i.e., $I \subset S$) to avoid issues regarding interpolation for discontinuous solutions. At each interface node, there exist two values of the solution, $\mathbf{u}_i^{-}$ and $\mathbf{u}_i^{+}$, representing the solution evaluated from the element of interest and the interface-adjacent element, respectively. The interface flux term can then be computed as 
\begin{equation}
    \mathbf{F}_{\partial \Omega_k} (\mathbf{u}_i) = \sum_{j\in I} \overline{\mathbf{F}}(\mathbf{u}_j^{-}, \mathbf{u}_j^{+}, \mathbf{n}_j) \overline{\phi}_j(\mathbf{x}),
\end{equation}
where $\overline{\mathbf{F}}(\mathbf{u}_i^{-}, \mathbf{u}_i^{+}, \mathbf{n}_i)$ are the common interface flux values dependent on the interior and exterior values of the solution and their associated normal vectors $\mathbf{n}_i$ and $\overline{\phi}_i(\mathbf{x})$ are the interface basis functions. The common interface flux is generally computed using an approximate Riemann solver such as that of \citet{Rusanov1962}. The interface basis functions are dependent on the choice of spatial discretization, e.g., for flux reconstruction schemes, these terms can be given as 
\begin{equation}
    \overline{\phi}_i(\mathbf{x}) = \mathbf{n}_i{\cdot}\mathbf{h}_i(\mathbf{x})- \phi_i(\mathbf{x}).
\end{equation}
Here, $\mathbf{h}_i$ are a set of correction functions \citep{Castonguay2011, Trojak2021} that posses the properties that 
\begin{equation}
    \mathbf{n}_i{\cdot} \mathbf{h}_j(\mathbf{x}_i) = \delta_{ij} \quad \mathrm{and} \quad \sum_{i \in I} \mathbf{h}_i (\mathbf{x}) \in \mathrm{RT}_{p},
\end{equation}
where $\mathrm{RT}_p$ is the Raviart--Thomas space \citep{Raviart1977} of order $p$. In this work, the flux reconstruction scheme with the equivalent discontinuous Galerkin correction functions \citep{Huynh2007} is used which recovers the nodal discontinuous Galerkin method \citep{Hesthaven2008DG}. The interface contribution to the divergence of the flux can then be given as
\begin{equation}\label{eq:cbarij}
    \boldsymbol{\nabla}{\cdot}{\mathbf{F}}_{\partial \Omega_k}(\mathbf{u}_i) = \sum_{j \in I} \overline{\mathbf{c}}_{ij} \overline{\mathbf{F}}(\mathbf{u}_j^{-}, \mathbf{u}_j^{+}, \mathbf{n}_j), \quad \mathrm{where} \quad \overline{\mathbf{c}}_{ij} = \nabla \overline{\phi}_i(\mathbf{x}_j).
\end{equation}
The semi-discrete form of \cref{eq:hyp_source} can then be given as
\begin{equation}
    \partial_t \mathbf{u}_i = -\left( \boldsymbol{\nabla}{\cdot}\mathbf{F}_{\partial \Omega_k} (\mathbf{u}_i) + \boldsymbol{\nabla}{\cdot}{\mathbf{F}}_{\partial \Omega_k}(\mathbf{u}_i)\right) + \mathbf{S_B}(\mathbf{u}_i).
\end{equation}

We assume that the spatial scheme satisfies the relation
\begin{equation}
    \partial_t \overline{\mathbf{u}} = - \int_{\partial \Omega_k}\overline{\mathbf{F}}\left(\mathbf{x} \right)\cdot\mathbf{n}(\mathbf{x})\ \mathrm{d}\mathbf{x} \approx - \sum_{j \in I} m_j \overline{\mathbf{F}}(\mathbf{u}_j^{-}, \mathbf{u}_j^{+}, \mathbf{n}_j)
\end{equation}
where $m_j$ is the associated quadrature weight for $\mathbf{x}_j$ and $\overline{\mathbf{u}}$ is the element-wise mean defined as
\begin{equation}
    \overline{\mathbf{u}} = \frac{1}{V_k}\int_{\Omega_k}\mathbf{u} (\mathbf{x})\ \mathrm{d}\mathbf{x} \quad \quad \mathrm{and} \quad \quad V_k = \int_{\Omega_k} \mathrm{d}\mathbf{x}.
\end{equation}
This assumption is appropriate for nodal discontinuous Galerkin schemes given appropriate quadrature and flux reconstruction schemes utilizing the equivalent discontinuous Galerkin correction functions.

\subsection{Eight-wave method}\label{ssec:powell}

A common method for enforcing a divergence-free magnetic field is to utilize the eight-wave method of \citet{Powell1999}. This approach relies on an additional wave structure of the Riemann problem in MHD that arises when the magnetic field is not exactly solenoidal, and it can be utilized to force the magnetic field to a solenoidal state via a source term, given as 
    \begin{equation}
        \mathbf{S}_{\boldsymbol{B}}(\mathbf{u}) = -\begin{bmatrix}
                0 \\ \boldsymbol{B} \\ \boldsymbol{u} \\ \boldsymbol{u}{\cdot}\boldsymbol{B}
            \end{bmatrix}\boldsymbol{\nabla}{\cdot}\boldsymbol{B}.
    \end{equation}
With the inclusion of this source term, the divergence of the magnetic field is typically suppressed to the order of magnitude of the approximation error \citep{Wu2018}. As such, due to the simplicity of implementation and applicability to general unstructured grids, it remains a routine approach for robustly enforcing the divergence-free constraint on the solenoidal field. In addition, only this modified form of the ideal MHD equations is symmetrizable and Galilean invariant when the magnetic field is not exactly solenoidal \citep{Wu2018}. However, as this form is non-conservative, it occasionally can cause inaccurate predictions around discontinuities in the flow (see \citet{Toth2000}). 

The use of Powell's method requires some clarification about the choice of the formulation for computing the divergence of the magnetic field. In the context of DSEM, there exist two formulations, a local divergence, consisting of just the interior component as
\begin{equation}
    \boldsymbol{\nabla}{\cdot}\boldsymbol{B}^L(\mathbf{u}_i) = \sum_{j \in S} \mathbf{c}_{ij}\boldsymbol{B}_j,
\end{equation}
and a global divergence, consisting of both the interior component and the interface contribution as
\begin{equation}
    \boldsymbol{\nabla}{\cdot}\boldsymbol{B}^G(\mathbf{u}_i) = \sum_{j \in S} \mathbf{c}_{ij}\boldsymbol{B}_j + \sum_{j \in I} \overline{\mathbf{c}}_{ij}\overline{\boldsymbol{B}}_j,
\end{equation}
where $\overline{\boldsymbol{B}}_j$ is a common interface value for the magnetic field, typically taken as the centered average of the interior and exterior values. Whereas the divergence-free constraint can be imposed on the local divergence through straightforward approaches such as projection to solenoidal bases, enforcing this constraint on the global divergence is typically more difficult as its domain of influence is not strictly contained within the element. It can be argued that the global approach is the ``correct'' choice as it is the one for which the space of the divergence is consistent with the space of the solution, but in practice, the local approach is typically sufficient. In this work, the global approach is used as the complexity of the two implementations is similar with Powell's method.
\section{Methodology}\label{sec:methodology}
Due to the presence of discontinuities in MHD flows in the form of hydrodynamic and magnetic shocks, it is necessary to apply some sort of a numerical stabilization procedure to ensure robustness of the DSEM approximation. In \citet{Dzanic2022}, an adaptive filtering approach was introduced with goal of stabilizing the scheme by discretely enforcing convex constraints on the solution, given in the form of 
\begin{equation}
    \Gamma(\mathbf{u}_i) > 0 \ \forall \ i \in S,
\end{equation}
where $\Gamma(\mathbf{u})$ is some constraint functional. For a positivity-preserving scheme, these constraints are set as 
\begin{equation}
    \Gamma_1(\mathbf{u}) = \rho \quad \mathrm{and} \quad \Gamma_2(\mathbf{u}) = P,
\end{equation}
corresponding to constraints on the positivity of density and pressure. 

While these constraints can ensure the positivity of these quantities, they are generally not restrictive enough to ensure that the solution remains well-behaved in the vicinity of discontinuities. It is necessary to attempt to form additional constraints on the solution that are restrictive enough to stabilize the solution in the vicinity of discontinuities without degrading the accuracy of the scheme in regions where the solution is smooth. By utilizing the fact that the minimum entropy principle presented in \cref{sec:preliminaries} should be satisfied by both smooth and discontinuous solutions, a third constraint on the solution is enforced corresponding to a discrete form of a local minimum entropy principle as
\begin{equation}
    \Gamma_3(\mathbf{u}) = \sigma(\mathbf{u}) - \sigma_{\min},
\end{equation}
where $\sigma(\mathbf{u}) = P\rho^{-\gamma}$ is the specific physical entropy and $\sigma_{\min}$ is some local minimum entropy bound. This minimum bound $\sigma_{\min}$ is computed in an element-wise manner as the discrete minima of the entropy functional across the element and its face neighbors prior to each stage of a temporal integration scheme, resulting in the enforcement of a discrete minimum entropy principle over the local domain of influence of the element (see \citet{Dzanic2022}, Section 2 and 3). It was found in the context of gas dynamics that enforcing this constraint ensured well-behaved solutions in the vicinity of discontinuities while recovering high-order accuracy in smooth regions of the flow \citep{Dzanic2022}.

\subsection{Adaptive filtering}
The constraints are enforced by an adaptive filtering procedure, where the filtered solution $\widetilde{\mathbf{u}}$ is given in terms of a filter kernel $H$ applied to the solution, i.e.,
\begin{equation}
    \widetilde{\mathbf{u}} = H(\mathbf{u}).
\end{equation}
This filtering is performed in modal space given a modal decomposition of the solution in the form of 
\begin{equation}
    \mathbf{u} (\mathbf{x}) = \sum_{i\in S} \widehat{\mathbf{u}}_i \psi_i(\mathbf{x}),
\end{equation}
where $\psi_i(\mathbf{x})\ \forall \ i \in S$ are a set of modal basis functions and $\widehat{\mathbf{u}}_i$ are their corresponding modes. We assume that this modal decomposition is chosen with respect to the unit measure (e.g., Legendre polynomials, Koornwinder polynomials, etc.). A discrete form of this change-of-basis operation can be given in terms of a Vandermonde matrix $\mathbf{V}$ as 
\begin{equation}
    \widehat{\mathbf{u}} = \mathbf{V}^{-1}\mathbf{u}.
\end{equation}
These change-of-basis operations and filtering procedures are applied in the reference coordinates to simplify the implementation.  

The filter kernel $\widehat{H}$ is taken as a second-order exponential kernel in modal space, such that the filtered modal modes can be computed as
\begin{equation}
    \widehat{H}_i(\widehat{\mathbf{u}}_i) = \widehat{\mathbf{u}}_i \exp (- \zeta p_i^2),
\end{equation}
where $\zeta$ is the filter strength and $p_i$ is the total order of the corresponding mode $\widehat{\mathbf{u}}_i$. It must be noted that the adaptive filtering approach is not restricted to this choice of filter and can be applied to any conservative filtering operation of one free variable that can recover both the unfiltered solution and the mean mode \citep{Dzanic2022}. The filtering operation $H(\mathbf{u})$ can be cast in terms of a matrix-vector operation as
\begin{equation}\label{eq:matvec}
    \widetilde{\mathbf{u}} = H(\mathbf{u}) = \mathbf{V} \boldsymbol{\Lambda} \mathbf{V}^{-1} \mathbf{u},
\end{equation}
where $\boldsymbol{\Lambda}$ is a diagonal matrix of $p+1$ unique values with its entries equal to $\boldsymbol{\Lambda}_{i,i} = \exp (- \zeta p_i^2)$.

The filter strength is computed via an element-wise nonlinear optimization process, taken as the minimum filter strength necessary such that the constraints are satisfied, i.e.,
\begin{equation}
     \zeta = \underset{\zeta\ \geq\ 0}{\mathrm{arg\ min}} \ \ \mathrm{s.t.} \ \  \left [\Gamma_1\left(\widetilde{\mathbf{u}}(\mathbf{x}_i) \right) > 0, \ \Gamma_2\left(\widetilde{\mathbf{u}}(\mathbf{x}_i) \right) > 0, \ \Gamma_3\left(\widetilde{\mathbf{u}}(\mathbf{x}_i) \right) > 0 \  \ \forall \ i \in S\right ].
\end{equation}
Existence of a solution of $\zeta$ is guaranteed if the element-wise mean of the solution satisfies the constraints, an assumption that will be explored in \cref{ssec:mhd}. As this optimization process is a function of a scalar free variable, its solution can be obtained using any root-bracketing approach. Furthermore, as it is nonlinear and non-convex, convergence to a local minima is sufficient in the case of multiple values of $\zeta$ existing such that the constraints are satisfied exactly. While this optimization problem seems computationally demanding due to the element-wise matrix-vector operations necessary to compute the filtered solution each iteration of the solve, we present a numerical approach to solving this problem in \cref{ssec:optimizations} that is much more computationally efficient than the original methodology in \citet{Dzanic2022}.

\subsection{Extensions to MHD}\label{ssec:mhd}
Extending the entropy filtering approach to the MHD system requires some modifications, with special care necessary in regards to the treatment of the source terms. The adaptive filtering operation naturally relies on that assumption that there exists a filter strength such that the constraints are satisfied, and it is trivial to show that a solution exists if the element-wise mean satisfies the constraints \citep{Dzanic2022}. The ability of discontinuous Galerkin-type approaches to preserve convex invariants of hyperbolic systems on the element-wise mean is a well established in the literature, and the reader is referred to a variety of works which utilize this property \citep{Zhang2010,Zhang2011,Zhang2011b,Zhang2011c,Chen2017,Wu2018,Dzanic2022}. However, the inclusion of the source term, the divergence-free condition, and the presence of entropy constraints introduces some caveats on this property of the scheme.

Let the set $G$ represent the set of solutions which satisfy the constraints (i.e., $\Gamma_1(\mathbf{u}) > 0, \Gamma_2(\mathbf{u}) > 0, \Gamma_3(\mathbf{u}) > 0$), and let the shorthand notation $\mathbf{u} \in G$ represent $\mathbf{u}_i \in G \ \forall \ i \in S$. To ensure that the filter can recover a constraint-satisfying solution, it is necessary for the temporal update of the element-wise mean to preserve these invariants, i.e., for some time step $n$, if $\mathbf{u}^n \in G$, then $\overline{\mathbf{u}}^{n+1} \in G$. For brevity, we consider a temporal update in the form of a forward Euler approximation, given as 
\begin{equation}
    \mathbf{u}^{n+1} = \mathbf{u}^n + \Delta t \left [L_1(\mathbf{u}^n) + L_2(\mathbf{u}^n) \right],
\end{equation}
where
\begin{equation}
    L_1(\mathbf{u}) = -\boldsymbol{\nabla}{\cdot}\mathbf{F}(\mathbf{u}) \quad \mathrm{and} \quad L_2(\mathbf{u}) = \mathbf{S_B}(\mathbf{u}).
\end{equation}
Without an exactly solenoidal magnetic field, the property $\overline{\mathbf{u}}^{n+1} \in G$ is not necessarily satisfied in this form under the standard assumptions posed in works such as \citet{Zhang2010} and the original presentation of entropy filtering for gas dynamics in \citet{Dzanic2022}, e.g., appropriate Riemann solver, CFL condition, strong stability preserving temporal integration. If we consider the set of solutions $G_P$ which satisfy just the positivity constraints (i.e.,  $\mathbf{u} \in G_P$ if $\Gamma_1(\mathbf{u}) > 0, \Gamma_2(\mathbf{u}) > 0$), then the work of \citet{Wu2018} showed that the property $\overline{\mathbf{u}}^{n+1} \in G_P$ is satisfied if the magnetic field is locally divergence-free. Furthermore, if we neglect the source term and consider an intermediate temporal update as
\begin{equation}
\mathbf{u}^{*} = \mathbf{u}^n + \Delta t L_1(\mathbf{u}^n),
\end{equation}
then the work of \citet{Bouchut2007} (paired with the equivalency of the element-wise mean and Godunov methods presented in \citet{Zhang2010} and subsequent works) shows that this intermediate state satisfies the property $\overline{\mathbf{u}}^{*} \in G$ under the standard assumptions. 

These two observations motivate an operator splitting approach for the filter. Two separate filtering operations are considered, a more restrictive filter which enforces both the positivity and entropy constraints, denoted by $H_e[\mathbf{u}]$, and a more relaxed filter that enforces only positivity constraints, denoted by $H_p[\mathbf{u}]$. As the assumption on the entropy constraints on the element-wise mean are satisfied by the intermediate state, the more restrictive filter can be applied, i.e.,
\begin{equation}
    \widetilde{\mathbf{u}}^{*} = H_e \left [\mathbf{u}^n + \Delta t L_1(\mathbf{u}^n)  \right ].
\end{equation}
Since the entropy constraints are the most restrictive constraint and the contribution of the source term is typically minimal compared to the divergence of the flux (since it is proportional to $\boldsymbol{\nabla}{\cdot} \mathbf{B}$), this filtering operation can usually mitigate the majority of the spurious oscillations in the vicinity of discontinuities. The contribution of the source term is then added onto this filtered state, after which the positivity constraints are then enforced again on the temporal update as
\begin{equation}
    \widetilde{\mathbf{u}}^{n+1} = H_p\left [\widetilde{\mathbf{u}}^{*} + \Delta t L_2(\mathbf{u}^n)  \right ].
\end{equation}
We remark here that it remains to be proven that element-wise mean for a nodal DG approximation preserves the positivity of the density and pressure if the magnetic field is not exactly locally divergence-free. As such, it is not guaranteed in a mathematical sense that there will exist a positivity-preserving filtered solution in the presence of small divergence errors. However, in the numerical experiments, the approach was found to be robust even in extreme flow conditions. 

Several properties of this splitting approach must be noted. First, it is very rarely the case that the secondary filtering operation is necessary -- the entropy constraints on $\widetilde{\mathbf{u}}^{*}$ are typically restrictive enough to where $\widetilde{\mathbf{u}}^{*} + \Delta t L_2(\mathbf{u}^n)$ retains its positivity-preserving properties, such that in most cases, the positivity constraints are typically just checked and no filtering is needed. Second, the splitting for the source term is calculated explicitly as $L_2 (\mathbf{u}^n)$, not through a Strang-type splitting approach \citep{Demkowicz1990} as $L_2 (\mathbf{u}^*)$. While the latter may potentially better approximate the necessary corrections to the solution for preserving a solenoidal magnetic field, these forms of splitting can introduce a limit on the temporal accuracy of the scheme and therefore are avoided. Finally, unless the linear filtering kernel which recovers the squeeze limiter of \citet{Zhang2010} is chosen (see \citet{Dzanic2022}, Remark 1), the divergence of the filtered magnetic field is not guaranteed to be equal or lower than the unfiltered state. As this work pertains to a nonlinear filter, it may introduce minor divergence errors similarly to any nonlinear limiting operation, but these are mitigated via the source term at the next temporal update with the explicit splitting approach such that its effects were found to be negligible. 

Extensions to higher-order strong stability preserving (SSP) schemes follow readily from this formulation, e.g., the temporal update for a third-order, three-stage SSP Runge--Kutta scheme, neglecting the notation $\widetilde{\cdot}$ for brevity, is given as
\begin{align}\label{eq:ssprk3}
\mathbf{u}^{(1)} &= H_p\left [H_e\left [\mathbf{u}^n + \Delta t L_1(\mathbf{u}^n) \right ]  + \Delta t L_2(\mathbf{u}^n)\right], \\
\mathbf{u}^{(2)} &= H_p\left [H_e\left [\frac{3}{4}\mathbf{u}^n + \frac{1}{4}\mathbf{u}^{(1)} + \frac{1}{4}\Delta t L_1(\mathbf{u}^{(1)}) \right ]  + \frac{1}{4}\Delta t L_2(\mathbf{u}^{(1)})\right], \nonumber \\
\mathbf{u}^{n+1} &= H_p\left [H_e\left [\frac{1}{3}\mathbf{u}^n + \frac{2}{3}\mathbf{u}^{(2)} + \frac{2}{3}\Delta t L_1(\mathbf{u}^{(2)}) \right ]  + \frac{2}{3}\Delta t L_2(\mathbf{u}^{(2)})\right], \nonumber
\end{align}
where the entropy constraints for $H_e$ are computed from the previous temporal stage (see \citet{Dzanic2022}, Appendix A).

\section{Implementation}\label{sec:implementation}

\begin{figure}[tbhp]
    \centering
    \adjustbox{width=0.33\linewidth, valign=b}{     \begin{tikzpicture}[spy using outlines={rectangle, height=3cm,width=2.3cm, magnification=3, connect spies}]
		\begin{axis}[name=plot1,
		    axis line style={draw=none},
		    tick style={draw=none},
		    axis x line=left,
            axis y line=left,
            axis equal image,
            clip mode=individual,
    		xmin=-1.6,
    		xmax=0,
    		xticklabels={,,},
    		ymin=-1,
    		ymax=1,
    		yticklabels={,,},
    		style={font=\Large},
    		scale = 1]
    		
            \draw [black, very thick] plot [] coordinates {(0, -1)  (0, 1) (-1.5, 0) (0, -1)};
            \fill[fill=black!10] (0, -1)--(0, 1)--(-1.5, 0);
            
            \draw [black, very thick, dotted] plot [] coordinates {(-.05, 1.1) (-1.55, 0.1) (-1.9, 1.6) (-.05, 1.1)};
            \draw [black, very thick, dotted] plot [] coordinates {(-1.58, -0.07) (0, -1.13) (-2.5, -1.13) (-1.58, -0.07)};
            \draw [black, very thick, dotted] plot [] coordinates {(-2.55, -1.0) (-2, 1.6) (-1.65, 0.03) (-2.55, -1.0)};
		  
            \draw[-,very thick, color=cyan!80!blue,fill=red] (0, -1) circle[radius=0.040];
            \draw[-,very thick, color=cyan!80!blue,fill=red] (0, 0) circle[radius=0.040];
            \draw[-,very thick, color=cyan!80!blue,fill=red] (0, 1) circle[radius=0.040];
            \draw[-,very thick, color=cyan!80!blue,fill=red] (-0.75, -0.5) circle[radius=0.040];
            \draw[-,very thick, color=cyan!80!blue,fill=red] (-1.5, 0) circle[radius=0.040];
            \draw[-,very thick, color=cyan!80!blue,fill=red] (-0.75, 0.5) circle[radius=0.040];
            \draw[-,very thick, color=black, fill=red] (-0.5, 0) circle[radius=0.040];
            
            \fill[-,fill=cyan!80!blue] (-.05, 1.1)  circle[radius=0.045];
            \fill[-,fill=cyan!80!blue] (-1.55, 0.1)  circle[radius=0.045];
            \fill[-,fill=cyan!80!blue] (-.8,0.6)  circle[radius=0.045];
            
            \fill[-,fill=cyan!80!blue] (-1.58, -0.07)  circle[radius=0.045];
            \fill[-,fill=cyan!80!blue] (-0.79, -0.6)  circle[radius=0.045];
            \fill[-,fill=cyan!80!blue] (0, -1.13) circle[radius=0.045];
            
            \node[] at (-0.2, -0.5) {$\Omega_k$};
            
        \end{axis} 		
	\end{tikzpicture}}
    \caption{\label{fig:scheme} Schematic of a two-dimensional $\mathbb P_2$ triangular element $\Omega_k$ showing interior solution points (red circles), interior interface flux/solution points (red circles, blue outline), and exterior interface flux points (blue circles).
    }
\end{figure}
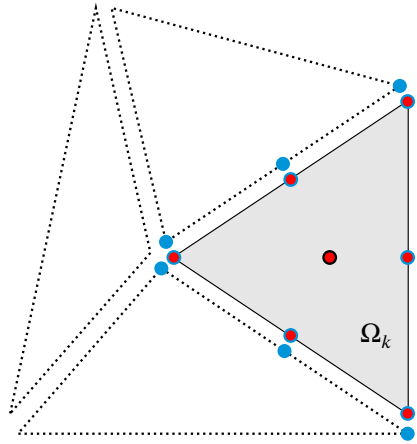

The governing equations and the adaptive filtering approach were implemented in PyFR \citep{Witherden2014}, a high-order GPU-accelerated unstructured flux reconstruction solver. The solution nodes were distributed along the Gauss--Legendre--Lobatto quadrature points and $\alpha$-optimized points \citep{Hesthaven2008DG} for tensor-product and simplex elements, respectively. An example of the solution and flux point distributions for a two-dimensional $\mathbb P_2$ triangular element is shown in \cref{fig:scheme}. Temporal integration was performed using a three-stage, third-order SSP Runge--Kutta scheme as presented in \cref{eq:ssprk3}. Unless otherwise stated, common interface fluxes were computed using the Harten-Lax-van Leer contact (HLLC) Riemann solver of \citet{Li2005b} and \citet{Gurski2004} with the Davis wavespeed estimate \citep{Davis1988}, although for most test cases, we observed negligible differences in comparison to Rusanov-type \citep{Rusanov1962} and Harten-Lax-van Leer (HLL) \citep{Harten1983} Riemann solvers. To avoid a vacuum state for the Riemann solver and apply a numerical tolerance to the entropy condition, the constraints were instead implemented as 
\begin{equation*}
    \Gamma_1(\mathbf{u}) = \rho - \epsilon, \quad \Gamma_2(\mathbf{u}) = P - \epsilon, \quad \mathrm{and} \quad  \Gamma_3(\mathbf{u}) = \sigma - \sigma_{\min} - \epsilon, 
\end{equation*}
where $\epsilon = 10^{-8}$ is a small constant taken as some arbitrary factor of the machine precision.

Boundary conditions were enforced in a weak sense through the imposition of an exterior ghost state to the interface Riemann solver \citep{Mengaldo2014}. Three types of boundary conditions were considered in this work: 1) Dirichlet boundary conditions, where the exterior state is explicitly defined; 2) Neumann boundary conditions, where the exterior state is identical to the interior state; and 3) reflecting boundary conditions, where the exterior state is identical to the interior state with the normal component of the velocity and magnetic field negated.

\subsection{Filter optimization}\label{ssec:optimizations}
Each time the filtering operation is called, the constraints are first checked on the solution. If the solution satisfies the constraints, no filtering is applied, otherwise the filter strength is computed using the Illinois root-bracketing approach \citep{Dowell1971} with a stopping tolerance of $10^{-8}$ and a maximum of 20 iterations. While the filter strength can be simply iterated by repeatedly evaluating the element-wise filtered solution as per \cref{eq:matvec} and computing the minima of the constraints, several optimizations can be performed to drastically decrease the computational cost of performing this filtering operation.  

First, instead of solving for $\zeta$, it beneficial to solve for $f = \exp (-\zeta)$ and utilize the relation
\begin{equation*}
    \exp(-\zeta p_i^2) = f^{p_i^2}.
\end{equation*}
This bounds the search space of the root-bracketing approach to $f \in [0,1]$, and the evaluation of the filter coefficients reduces to simple integer powers of the argument $f$. Then, to avoid the costly computation of the matrix-vector product in \cref{eq:matvec} each iteration of the root-bracketing process, certain properties of the matrix $\boldsymbol{\Lambda}$ can be exploited. As previously mentioned, $\boldsymbol{\Lambda}$ is a diagonal matrix of $p+1$ unique values with its entries equal to $\boldsymbol{\Lambda}_{i,i} = \exp (- \zeta p_i^2)$. If we define a set of diagonal matrices $\mathbf{I}^{(k)}$ for $0 \leq k \leq p$ as
\begin{equation}
    \mathbf{I}^{(k)}_{i,i} = \begin{cases}
    1, \quad \mathrm{if} p_i = k, \\
    0, \quad \mathrm{else},
    \end{cases}
\end{equation}
then the filtering operation can be equivalently represented as
\begin{equation}
    \widetilde{\mathbf{u}} = \sum_{i=0}^p f^{p_i^2} \mathbf{u}^{(k)},
\end{equation}
where
\begin{equation}
    \mathbf{u}^{(k)} = \mathbf{V} \mathbf{I}^{(k)} \mathbf{V}^{-1} \mathbf{u}.
\end{equation}
Note that the values $\mathbf{u}^{(k)}$ are independent of the value of $f$, such that these values can be pre-computed and the filtered solution can be efficiently evaluated each iteration of the root-bracketing approach without having to repeatedly compute the matrix-vector product $\mathbf{V} \boldsymbol{\Lambda} \mathbf{V}^{-1} \mathbf{u}$. 

This approach can be even further optimized by utilizing the fact that the nodal values of the solution can now be decoupled, such that the root-bracketing process can be applied across each solution node sequentially which is particularly beneficial for computing architectures where memory bandwidth is the bottleneck. In this sequential approach, each solution node $\mathbf{x}_j$ for $j \in S$ solves for a value of $f_j$ such that $\widetilde{\mathbf{u}}_j$ satisfies the constraints. It is trivial to show that if 
\begin{equation*}
    f = \underset{j \in S}{\min}\ f_j,
\end{equation*}
then $\widetilde{\mathbf{u}}$ satisfies the constraints at all nodes. It is therefore advantageous to then use $f_j$ as the upper bound for the root-bracketing process for the node $\mathbf{x}_{j+1}$ as the constraints can be checked for $\widetilde{\mathbf{u}}_{j+1}$ using this upper bound and the root-bracketing process for that node can be skipped if they are satisfied. As the proposed algorithm requires effectively only one full evaluation of \cref{eq:matvec} irrespective of the number of iterations of the root-bracketing approach, the memory bandwidth requirements are significantly decreased, such that the filtering process becomes only a relatively small portion of the total compute time that is typically less than the cost of the evaluation of the divergence of the flux. An example of the implementation of this approach is provided in the electronic supplementary material of this work, and an evaluation of the efficiency improvements of this proposed algorithm in comparison to the original methodology in \citet{Dzanic2022} which utilizes repeated evaluations of \cref{eq:matvec} is presented in \cref{sec:results}. 

\section{Results}\label{sec:results}
\subsection{Near-vacuum convecting vortex}
To verify that the proposed scheme retains the high-order accuracy of the underlying DSEM for smooth solutions, the rate of convergence was calculated for the smooth magnetized convecting vortex problem introduced by \citet{Christlieb2015}. For this problem, the domain is taken as $\Omega = [-10, 10]^2$ with periodic boundary conditions discretized on a structured quadrilateral mesh, and the initial conditions are given as
\begin{equation}
    \mathbf{q}(\mathbf{x}, 0) =    
    \begin{bmatrix}
        \rho \\
        u \\
        v \\ 
        B_x \\
        B_y \\
        P
    \end{bmatrix} =
    \begin{bmatrix}
        1 \\
        1 - y \delta u \\
        1 + x \delta u \\
        -y \delta B \\
        \phantom{-}x \delta B \\
        1 + \delta P
    \end{bmatrix}
\end{equation}
where
\begin{equation}
    \delta u = \frac{\mu}{\sqrt{2}\pi} \phi(r), \quad \quad \delta B = \frac{\mu}{2\pi} \phi(r), \quad \quad \delta P = -\frac{\mu^2 (1 + r^2)}{8 \pi^2}\phi(r)^2,
\end{equation}
and
\begin{equation}
  \phi(r) = \exp (1 - r^2),  \quad \quad r = \sqrt{x^2 + y^2}.
\end{equation}
The specific heat ratio was set as $\gamma = 5/3$. These conditions give a non-isentropic nature to the flow field, which allows for a proper assessment of the proposed entropy-based constraints for smooth flows where the filter should be primarily inactive. To make this problem more numerically challenging, the parameter $\mu$ is chosen such as to give a near-vacuum state for the pressure field \citep{Christlieb2015}. This value was set as $\mu = 5.38948938512$, which gives a minimum pressure value in the domain of approximately $2 \epsilon = 2{\cdot}10^{-8}$.

The problem was solved until a non-dimensional time of $t = 0.05$ using a fixed time step of $\Delta t = 1{\cdot}10^{-4}$, after which the $L^1$ norm of the magnetic field error was computed as
\begin{equation}
    e_B = \frac{1}{A} \int_{\Omega} \left | B_x- B^{\mathrm{exact}}_x \right |  + \left | B_y- B^{\mathrm{exact}}_y \right |  \ \mathrm{d}\mathbf{x},
\end{equation}
where $A = 20^2$. The exact solution was computed through a translation of the initial conditions with a translation velocity of $[1, 1]$, and the integration was computed using a 9th order Gauss--Legendre quadrature rule. The error with respect to the mesh resolution $N_e$ is shown for various approximation orders in \cref{tab:vortex_error} in addition to the rate of convergence. High-order convergence, on the order of $p$ to $p+1$, was observed for all approximation orders between $\mathbb P_2$ and $\mathbb P_5$. Furthermore, the $\mathbb P_2$ results can be compared to the positivity-preserving third-order DG scheme in \citet{Wu2018} (Table 2), which can be considered as a subset of the adaptive filtering approach without entropy constraints (see \citet{Dzanic2022}, Remark 1). The proposed scheme gives marginally lower error even after removing the $1/A$ normalization factor. The maximum value of the limiting factor $1 - f$ across the domain over the simulation time is also shown in \cref{tab:vortex_f} for the various approximation orders and mesh resolution. It can be seen that the effect of the filtering is essentially negligible for smooth solutions, with values on the order of $10^{-4}$ to $10^{-7}$.

\begin{figure}[tbhp]
    \centering
    \begin{tabularx}{0.8\textwidth}{r  @{\extracolsep{\fill}} cccc}
        $N_e$ & $\mathbb{P}_2$ &$\mathbb{P}_3$ &$\mathbb{P}_4$ &$\mathbb{P}_5$ \\ 
        \midrule
        $20^2$ & \num{2.29e-04} & \num{3.07e-05} & \num{3.30e-06} & \num{4.04e-07} \\
        $25^2$ & \num{1.18e-04} & \num{1.17e-05} & \num{1.37e-06} & \num{9.36e-08} \\
        $33^2$ & \num{5.39e-05} & \num{4.07e-06} & \num{3.39e-07} & \num{2.44e-08} \\
        $40^2$ & \num{3.04e-05} & \num{2.02e-06} & \num{1.44e-07} & \num{9.55e-09} \\
        $50^2$ & \num{1.59e-05} & \num{8.72e-07} & \num{5.02e-08} & \num{2.92e-09} \\
        $67^2$ & \num{6.92e-06} & \num{2.99e-07} & \num{1.31e-08} & \num{6.51e-10} \\
        \midrule
        \textbf {RoC} & $2.89$ & $3.80$ & $4.62$ & $5.23$
    \end{tabularx}
    \captionof{table}{\label{tab:vortex_error} Convergence of the $L^1$ norm of the magnetic field error at $t = 0.05$ with respect to mesh resolution $N_e$ for the near-vacuum convecting vortex problem with varying approximation order. Rate of convergence shown on bottom.}
\end{figure}

\begin{figure}[tbhp]
    \centering
    \begin{tabularx}{0.8\textwidth}{r  @{\extracolsep{\fill}} cccc}
        $N_e$ & $\mathbb{P}_2$ &$\mathbb{P}_3$ &$\mathbb{P}_4$ &$\mathbb{P}_5$ \\ 
        \midrule
        $20^2$ & \num{1.08e-04} & \num{1.64e-05} & \num{1.89e-05} & \num{5.38e-07} \\
        $25^2$ & \num{3.30e-08} & \num{9.99e-06} & \num{3.79e-08} & \num{1.09e-06} \\
        $33^2$ & \num{5.92e-08} & \num{7.85e-06} & \num{6.41e-08} & \num{9.66e-07} \\
        $40^2$ & \num{7.69e-05} & \num{2.09e-05} & \num{2.46e-06} & \num{3.38e-07} \\
        $50^2$ & \num{6.39e-05} & \num{1.74e-05} & \num{1.34e-06} & \num{2.93e-07} \\
        $67^2$ & \num{2.51e-07} & \num{4.21e-06} & \num{2.55e-07} & \num{7.53e-07} \\
        \midrule
    \end{tabularx}
    \captionof{table}{\label{tab:vortex_f} Maximum value of the limiting factor $1 - f$ in the domain over the simulation time with respect to mesh resolution $N_e$ for the near-vacuum convecting vortex problem with varying approximation order.}
\end{figure}

\subsection{Brio--Wu shock tube}

Extensions to flows with discontinuities was then performed through the shock tube problem of \citet{Brio1988} which includes features of the Riemann problem such as shock waves, contact discontinuities, rarefaction waves, and compound waves. For this problem, the domain is set as $\Omega = [0,1]$ and the initial conditions are given by

\begin{equation}
    \mathbf{q}(\mathbf{x}, 0) =
    \begin{bmatrix}
        \rho \\
        u \\
        v \\ 
        B_x \\
        B_y \\
        P
    \end{bmatrix}
    = \begin{cases}
    \mathbf{q}_l, \quad \mathrm{if} \ x \leq 0.5, \\
    \mathbf{q}_r, \quad \mathrm{else},
    \end{cases}
    \quad \mathrm{where} \quad 
    \mathbf{q}_l
    =
    \begin{bmatrix}
        1 \\
        0 \\
        0 \\
        0.75 \\
        1 \\
        1
    \end{bmatrix}
    \quad \mathrm{and} \quad \mathbf{q}_r
    =
    \begin{bmatrix}
        0.125 \\
        0 \\
        0 \\
        0.75 \\
        -1 \\
        0.1
    \end{bmatrix}.
\end{equation}
The specific heat ratio is set as $\gamma = 2$. The hydrodynamic components of this problem are identical to the Sod shock tube \citep{Sod1978}. Although this problem is one-dimensional, it was instead solved on a one element wide two-dimensional quadrilateral mesh to facilitate the use of the vertical magnetic field component within the solver. Dirichlet boundary conditions were applied on the left/right boundaries while periodic boundary conditions were applied along the top/bottom boundaries. 

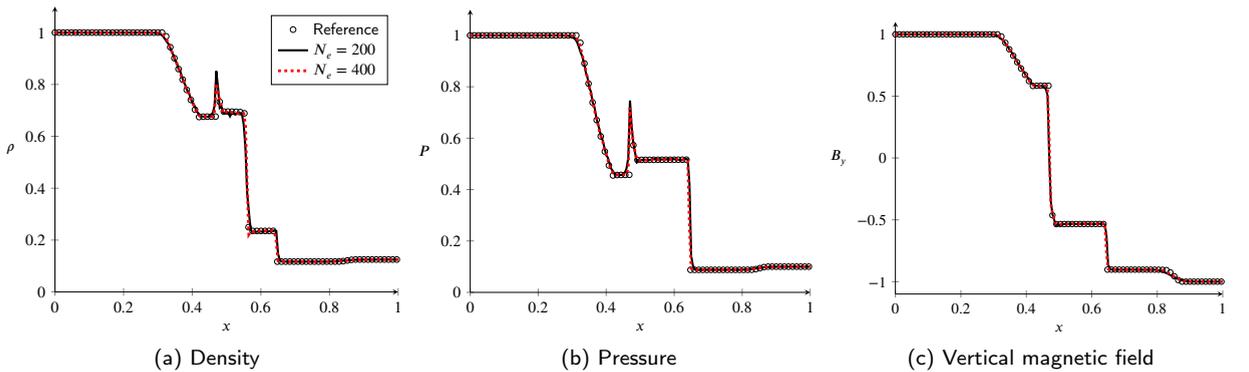
\begin{figure}[tbhp]
    \centering
    \subfloat[Density]{\adjustbox{width=0.33\linewidth, valign=b}{\begin{tikzpicture}[spy using outlines={rectangle, height=3cm,width=2.3cm, magnification=3, connect spies}]
	\begin{axis}[name=plot1,
		axis line style={latex-latex},
	    axis x line=left,
        axis y line=left,
		xlabel={$x$},
    	xmin=0, xmax=1,
    	xtick={0, 0.2, 0.4, 0.6, 0.8, 1.0},
    	ylabel={$\rho$},
    	ymin=0,ymax=1.1,
    	ytick={0, 0.2, 0.4, 0.6, 0.8, 1.0},
        clip mode=individual,
    	ylabel style={rotate=-90},
    	legend style={at={(0.97, 0.97)},anchor=north east},
    	legend cell align={left},
     ]
    	
        \addplot[color=black, style={ultra thin}, only marks, mark=o, mark options={scale=0.8}, mark repeat = 3, mark phase =0]
        table[x expr={\thisrow{x}+0.5}, y = rho, col sep=comma]{./figs/data/briowu_shocktube_ref_250.csv};
        \addlegendentry{Reference};

        \addplot[color=black, 
                style={very thick}]
        table[x expr={\thisrow{x}+0.5}, y=rho, col sep=comma]{./figs/data/bw200.csv};
        \addlegendentry{$N_e = 200$};
        
        \addplot[color=red, 
                style={ultra thick, dotted}]
        table[x expr={\thisrow{x}+0.5}, y=rho, col sep=comma]{./figs/data/bw400.csv};
        \addlegendentry{$N_e = 400$} ;
	\end{axis}
\end{tikzpicture}}}
    \subfloat[Pressure]{\adjustbox{width=0.33\linewidth, valign=b}{\begin{tikzpicture}[spy using outlines={rectangle, height=3cm,width=2.3cm, magnification=3, connect spies}]
	\begin{axis}[name=plot1,
		axis line style={latex-latex},
	    axis x line=left,
        axis y line=left,
		xlabel={$x$},
    	xmin=0, xmax=1,
    	xtick={0, 0.2, 0.4, 0.6, 0.8, 1.0},
    	ylabel={$P$},
    	ymin=0,ymax=1.1,
    	ytick={0, 0.2, 0.4, 0.6, 0.8, 1.0},
        clip mode=individual,
    	ylabel style={rotate=-90},
    	legend style={at={(0.97, 0.97)},anchor=north east},
    	legend cell align={left},
     ]
    	
        \addplot[color=black, style={ultra thin}, only marks, mark=o, mark options={scale=0.8}, mark repeat = 3, mark phase =0]
        table[x expr={\thisrow{x}+0.5}, y = p, col sep=comma]{./figs/data/briowu_shocktube_ref_250.csv};
        
        \addplot[color=black, 
                style={very thick}]
        table[x expr={\thisrow{x}+0.5}, y=p, col sep=comma]{./figs/data/bw200.csv};
        
        \addplot[color=red, 
                style={ultra thick, dotted}]
        table[x expr={\thisrow{x}+0.5}, y=p, col sep=comma]{./figs/data/bw400.csv} ;

	\end{axis}
\end{tikzpicture}}}
    \subfloat[Vertical magnetic field]{\adjustbox{width=0.33\linewidth, valign=b}{\begin{tikzpicture}[spy using outlines={rectangle, height=3cm,width=2.3cm, magnification=3, connect spies}]
	\begin{axis}[name=plot1,
		axis line style={latex-latex},
	    axis x line=left,
        axis y line=left,
		xlabel={$x$},
    	xmin=0, xmax=1,
    	xtick={0, 0.2, 0.4, 0.6, 0.8, 1.0},
    	ylabel={$B_y$},
    	ymin=-1.1,ymax=1.1,
    	ytick={-1, -.5, 0, 0.5, 1},
        clip mode=individual,
    	ylabel style={rotate=-90},
    	legend style={at={(0.97, 0.97)},anchor=north east},
    	legend cell align={left},
     ]
    	
        \addplot[color=black, style={ultra thin}, only marks, mark=o, mark options={scale=0.8}, mark repeat = 3, mark phase =0]
        table[x expr={\thisrow{x}+0.5}, y = By, col sep=comma]{./figs/data/briowu_shocktube_ref_250.csv};
        
        \addplot[color=black, 
                style={very thick}]
        table[x expr={\thisrow{x}+0.5}, y=By, col sep=comma]{./figs/data/bw200.csv};
        
        \addplot[color=red, 
                style={ultra thick, dotted}]
        table[x expr={\thisrow{x}+0.5}, y=By, col sep=comma]{./figs/data/bw400.csv} ;

	\end{axis}
\end{tikzpicture}}}
    \newline
    \caption{
    \label{fig:bwshock} Density, pressure, and vertical magnetic field profiles for the Brio--Wu shock tube problem at $t = 0.1$ computed using a $\mathbb P_3$ FR scheme with $200$ and $400$ elements. 
    }
\end{figure}

The problem was computed with a $\mathbb P_3$ scheme using a coarser mesh of 200 elements and a finer mesh of 400 elements with time steps of $\Delta t = 2{\cdot}10^{-4}$ and $1{\cdot}10^{-4}$, respectively. A reference solution was also computed using a highly-resolved $\mathbb P_0$ scheme with $5{\cdot}10^4$ elements. The predicted density, pressure, and vertical magnetic field profiles at $t = 0.1$ are shown in \cref{fig:bwshock} for both the coarse and fine mesh. For all fields, both rarefaction waves and the shock were well-resolved, showing sub-element resolution without any noticeable spurious oscillations. Furthermore, similar behavior was observed for the contact and compound wave in the pressure and magnetic fields. Some minor oscillations were observed in the density profile in the region between the compound wave and contact discontinuity, although this behavior is not uncommon for some numerical schemes. The predicted density profile in that region converged to the reference results with increasing resolution, but minor undershoots in front of the contact discontinuity were observed. 

\subsection{Orszag--Tang vortex}
Two-dimensional flows with more complex features were then considered through the canonical Orszag--Tang vortex problem \citep{Orszag1979}. This case is a well-known model problem for evaluating a scheme's ability to handle MHD shocks and shock interactions as well as predicting transition to supersonic MHD turbulence. The domain is set as $\Omega = [0,1]^2$ with periodic boundary conditions, and the initial conditions are given as    
\begin{equation}
    \mathbf{q}(\mathbf{x}, 0) =     
    \begin{bmatrix}
        \rho \\
        u \\
        v \\ 
        B_x \\
        B_y \\
        P
    \end{bmatrix}
    =
    \begin{bmatrix}
        25/(36\pi) \\
        -\sin(2\pi y) \\
        \phantom{-}\sin(2\pi x) \\
        \phantom{-}\sin (2 \pi y)/\sqrt{4 \pi}\\
        -\sin (4 \pi x)/\sqrt{4 \pi}\\
        5/(12\pi)
    \end{bmatrix}.
\end{equation}

\begin{figure}[htbp!]
        \centering
        \subfloat[$N_e = 64^2$]{\label{fig:OTV_64} \adjustbox{width=0.389\linewidth,valign=b}{\includegraphics[width=\textwidth]{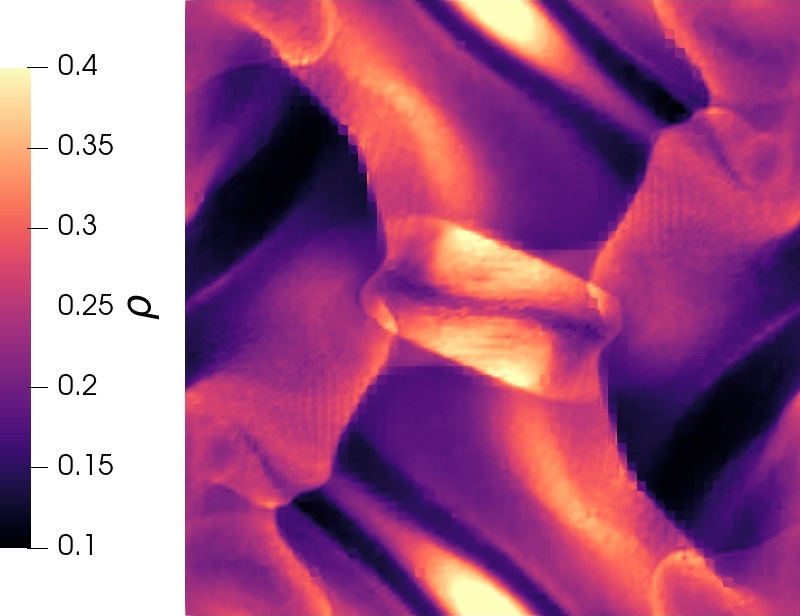}}}
        ~
        \subfloat[$N_e = 128^2$]{\label{fig:OTV_128} \adjustbox{width=0.3\linewidth,valign=b}{\includegraphics[width=\textwidth]{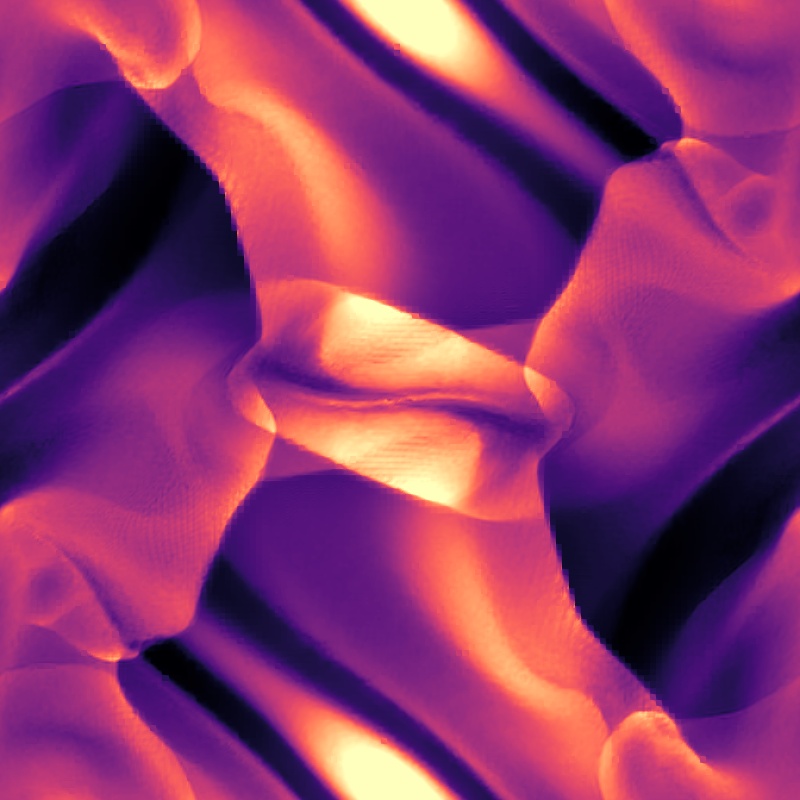}}}
        ~
        \subfloat[$N_e = 256^2$]{\label{fig:OTV_256} \adjustbox{width=0.3\linewidth,valign=b}{\includegraphics[width=\textwidth]{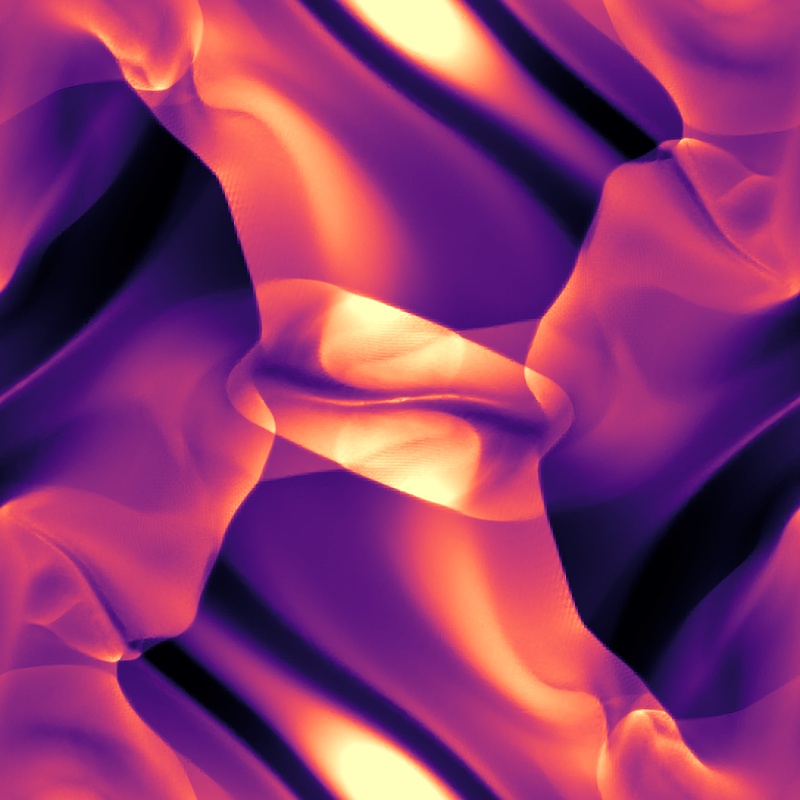}}}
        \newline
        \caption{\label{fig:OTV_d} Contours of density for the Orszag-Tang vortex at $t = 0.5$ computed using a $\mathbb P_3$ FR scheme with $64^2$ (left), $128^2$ (middle), and $256^2$ (right) elements.}
    \end{figure}

\begin{figure}[htbp!]
        \centering
        \subfloat[$N_e = 64^2$]{\adjustbox{width=0.389\linewidth,valign=b}{\includegraphics[width=\textwidth]{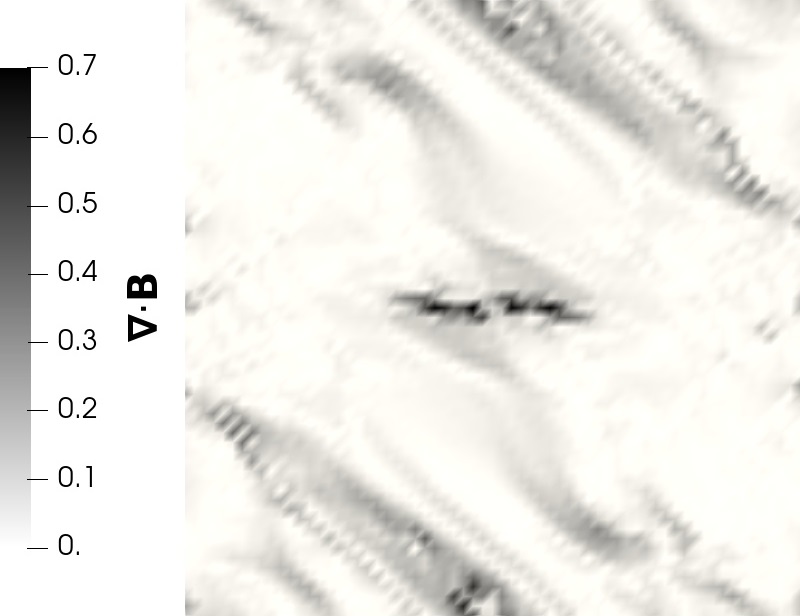}}}
        ~
        \subfloat[$N_e = 128^2$]{\adjustbox{width=0.3\linewidth,valign=b}{\includegraphics[width=\textwidth]{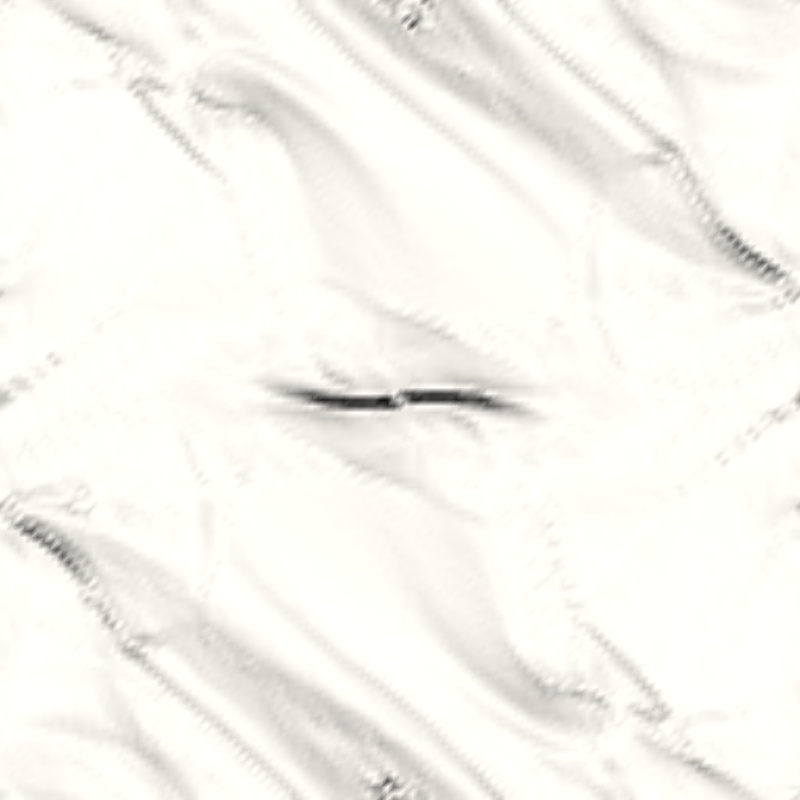}}}
        ~
        \subfloat[$N_e = 256^2$]{\adjustbox{width=0.3\linewidth,valign=b}{\includegraphics[width=\textwidth]{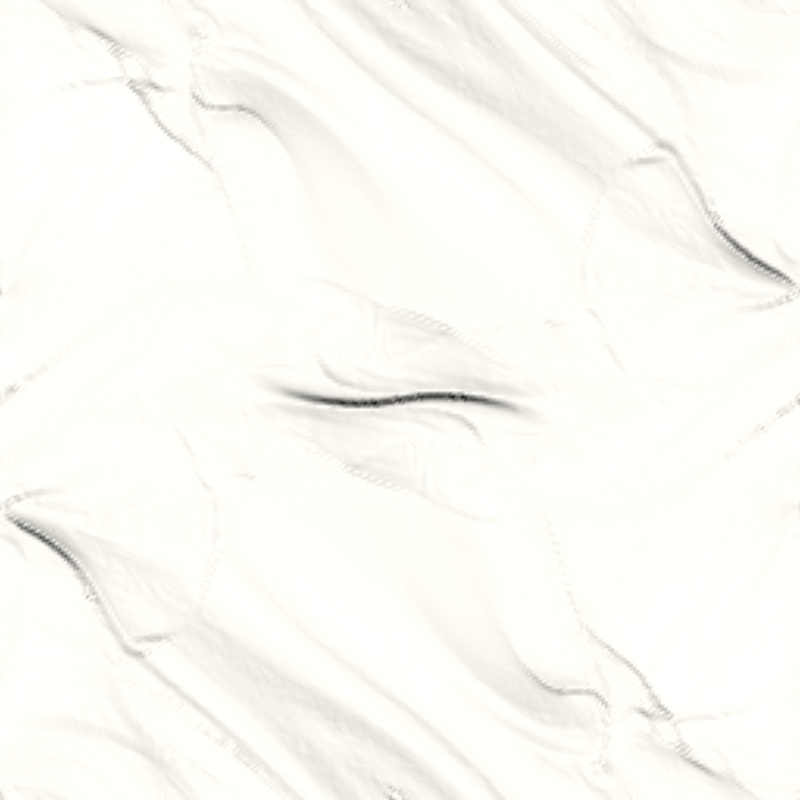}}}
        \newline
        \caption{\label{fig:OTV_divb} Contours of absolute magnetic divergence for the Orszag-Tang vortex at $t = 0.5$ computed using a $\mathbb P_3$ FR scheme with $64^2$ (left), $128^2$ (middle), and $256^2$ (right) elements.}
    \end{figure}

\begin{figure}[htbp!]
        \centering
        \subfloat[$N_e = 64^2$]{\adjustbox{width=0.389\linewidth,valign=b}{\includegraphics[width=\textwidth]{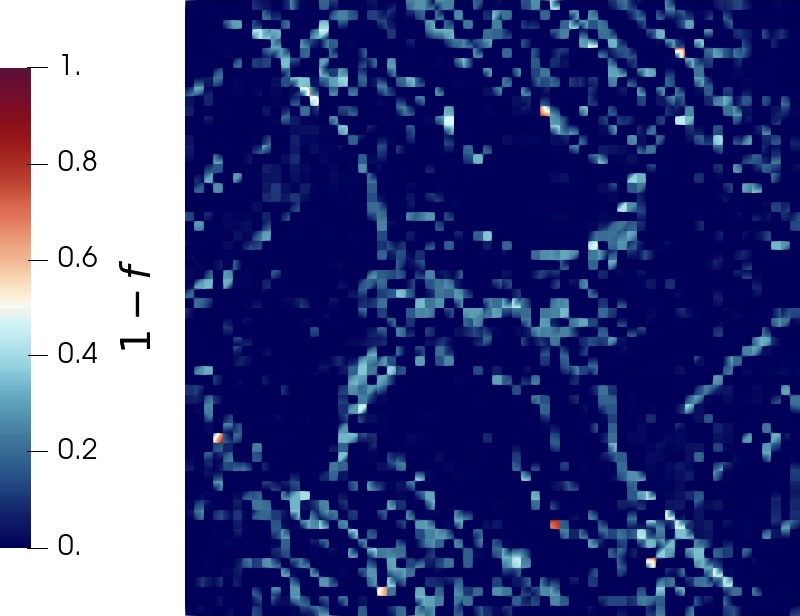}}}
        ~
        \subfloat[$N_e = 128^2$]{\adjustbox{width=0.3\linewidth,valign=b}{\includegraphics[width=\textwidth]{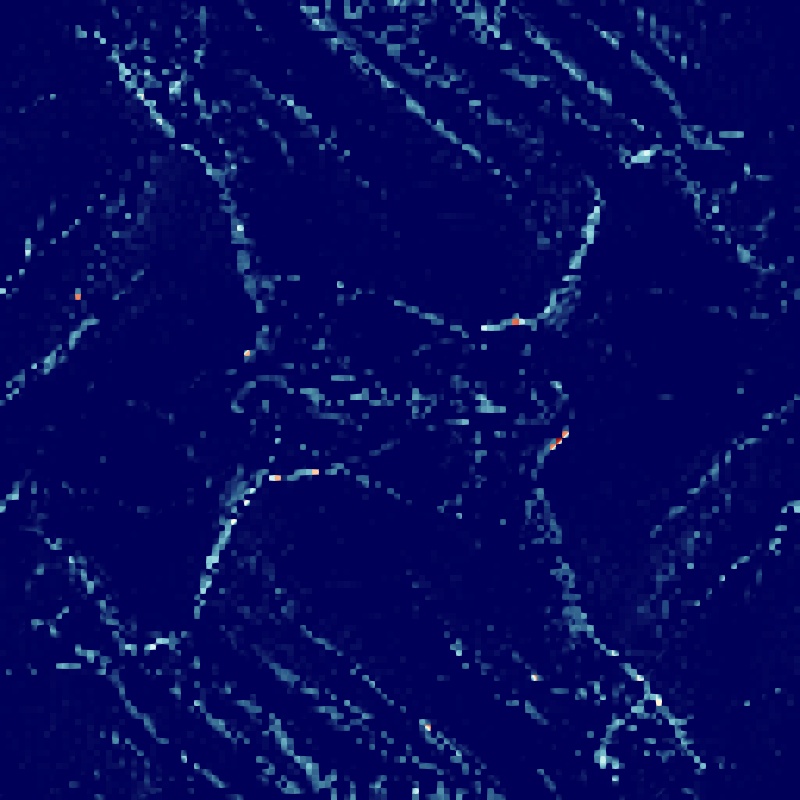}}}
        ~
        \subfloat[$N_e = 256^2$]{\adjustbox{width=0.3\linewidth,valign=b}{\includegraphics[width=\textwidth]{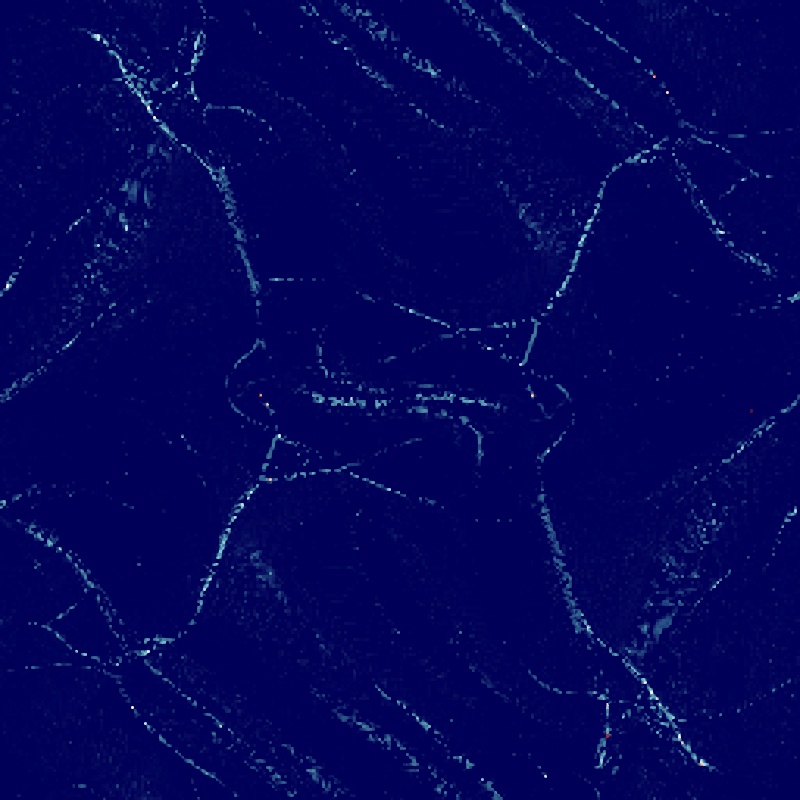}}}
        \newline
        \caption{\label{fig:OTV_f} Contours of the limiting factor $1 - f$ for the Orszag-Tang vortex at $t = 0.5$ computed using a $\mathbb P_3$ FR scheme with $64^2$ (left), $128^2$ (middle), and $256^2$ (right) elements.}
    \end{figure}

The specific heat ratio is set as $\gamma = 5/3$. Uniform quadrilateral meshes of various resolution were generated, and the problem was solved using a $\mathbb P_3$ scheme. The contours of density at $t = 0.5$ computed on meshes with $N_e = 64^2$, $128^2$, and $256^2$ elements are shown in \cref{fig:OTV_d}, computed using time steps of $\Delta t = 4{\cdot}10^{-4}$, $2{\cdot}10^{-4}$, and $1{\cdot}10^{-4}$, respectively. The results show good prediction of the canonical flow field of the Orszag--Tang vortex, with better approximation of shock structure and small-scale flow features with increasing resolution. Minor spurious oscillations were observed in the density field at low resolutions, but these oscillations diminished with increasing mesh resolution, such that the flow field at $N_e = 256^2$ was virtually oscillation-free. Furthermore, the contours of the absolute divergence of the magnetic field are shown in \cref{fig:OTV_divb} for the three meshes. It can be seen that the divergence errors tend to be focused in areas with large flow gradients and that the magnitude of the divergence errors decreases with increasing mesh resolution. Additionally, the distribution of the limiting factor $1 - f$ is shown in \cref{fig:OTV_f}. Limiting was typically focused around regions with discontinuities, with some minor activation of the filter in smooth regions. As the resolution was increased, these spurious activations decreased, such that the filter became increasingly focused in the discontinuous regions. 

\begin{figure}[htbp!]
        \centering
        \adjustbox{width=0.8\linewidth,valign=b}{    \begin{tikzpicture}[spy using outlines={rectangle, height=3cm,width=2.3cm, magnification=3, connect spies}]
		\begin{axis}[name=plot1,
            width=\textwidth,
            height=0.8*\axisdefaultheight,
		    axis line style={latex-latex},
		    axis x line=left,
            axis y line=left,
            clip mode=individual,
		    xlabel={$x/L$},
		    xtick={0, 0.2, 0.4, 0.6, 0.8, 1},
    		xmin=0,
    		xmax=1,
    	    x tick label style={
        		/pgf/number format/.cd,
            	fixed,
            	fixed zerofill,
            	precision=1,
        	    /tikz/.cd},
    		ylabel={$p$},
    		ylabel style={rotate=-90},
    		ytick={0, 0.1, 0.2,0.3},
    		ymin=0,
    		ymax=0.3,
    		y tick label style={
        		/pgf/number format/.cd,
            	fixed,
            	fixed zerofill,
            	precision=1,
        	    /tikz/.cd},
    		legend style={at={(0.03,0.03)},anchor=south west,font=\small},
    		legend cell align={left},
    		style={font=\normalsize}]

            \addplot[color=black, style={ultra thin}, only marks, mark=o, mark options={scale=0.8}, mark repeat = 2, mark phase =0]
				table[x index=0,y index=1,col sep=comma,unbounded coords=jump]{./figs/data/otv_ref_p_jiang_wu_1999.csv};
			\addlegendentry{ Reference}
				
			\addplot[color=black, style={ultra thick, dotted}]
				table[x=x,y=p,col sep=comma,unbounded coords=jump]{./figs/data/otv64_slice.csv};
			\addlegendentry{$N_e = 64^2$}
   
			\addplot[color=black, style={very thick}]
				table[x=x,y=p,col sep=comma,unbounded coords=jump]{./figs/data/otv128_slice.csv};
			\addlegendentry{$N_e = 128^2$}
   
			\addplot[color=red!90!black, style={very thick}]
				table[x=x,y=p,col sep=comma,unbounded coords=jump]{./figs/data/otv256_slice.csv};
			\addlegendentry{$N_e = 256^2$ }

		\end{axis} 		
	\end{tikzpicture}}
        \caption{\label{fig:OTV_p} 
        Pressure profile of the Orszag-Tang vortex on the cross-section $y/L = 0.3125$ at $t = 0.48$ computed using a $\mathbb P_3$ FR scheme with $64^2$, $128^2$, and $256^2$ elements. Numerical results of \citet{Jiang1999} shown for reference.
        }
    \end{figure}
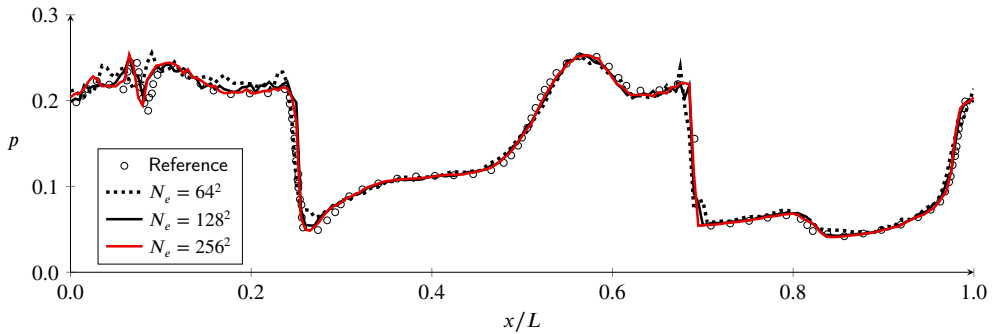

For a more quantitative comparison, the predicted pressure profile on the cross-section $y = 0.3125$ at $t = 0.48$ is shown in \cref{fig:OTV_p} in comparison to the results of \citet{Jiang1999} obtained using a high-order weighted essentially non-oscillator (WENO) scheme. It can be seen that similar observations can be made for the pressure field as with the density field, with minor spurious oscillations at lower resolutions that diminish with increasing resolution. The overall prediction of the pressure profile was in good agreement with the reference results at moderate and high resolutions, and strong discontinuities in the pressure field were generally well resolved at the sub-element level even at lower resolutions. At the highest resolution, there was very good agreement with the reference data with minor differences in the prediction of the location of some small-scale flow features on the left-hand side of the cross-section. Overall, the results showed good predictions of the various flow features of the Orszag--Tang vortex. 

\subsection{Shock cloud interaction}
The proposed scheme was then evaluated on the shock cloud interaction problem of \citet{Dai1998}, consisting of a high-density cloud interacting with an impinging shock wave which results in strong discontinuities and the development of small-scale flow instabilities. The problem setup as described by \citet{Balbs2006} is solved on the domain $\Omega = [0, 1]^2$ with the initial conditions given as
\begin{equation}
    \mathbf{q}(\mathbf{x}, 0) =  
       [\rho, 
        u, 
        v, 
        w, 
        B_x, 
        B_y,
        B_z,
        P]^T
    = \begin{cases}
    \mathbf{q}_c, \quad \mathrm{if} \ r' \leq 0.15, \\
    \mathbf{q}_l, \quad \mathrm{else\ if} \ x \leq 0.6, \\
    \mathbf{q}_r, \quad \mathrm{else},
    \end{cases}
\end{equation}
where the cloud state, left state, and right state are given as
\begin{equation}
    \mathbf{q}_c
    =
    \begin{bmatrix}
        10 \\
        0 \\
        0 \\
        0 \\
        \phantom{-}2.1826182 \\
        -2.1826182 \\
        1 \\
        167.345
    \end{bmatrix}, 
    \quad
    \mathbf{q}_l
    =
    \begin{bmatrix}
        3.86859 \\
        0 \\
        0 \\
        0 \\
        \phantom{-}2.1826182 \\
        -2.1826182 \\
        1 \\
        167.345
    \end{bmatrix}, 
    \quad \mathrm{and} \quad \mathbf{q}_r
    =
    \begin{bmatrix}
        1 \\
        -11.2536 \\
        0 \\
        0 \\
        0 \\
        0.56418958 \\
        0.56418958 \\
        1
    \end{bmatrix},
\end{equation}
respectively. The specific heat ratio is set as $\gamma = 5/3$. The cloud is centered at $[0.8, 0.5]$ with a radius of $0.15$, such that 
\begin{equation*}
r' = \sqrt{(x - 0.8)^2 + (y - 0.5)^2}.
\end{equation*}

To facilitate the use of the three-dimensional magnetic field within the solver, the problem is solved on a one element deep three-dimensional hexahedral mesh. Additionally, while the original problem setup uses a $[0, 1]^2$ domain with Neumann boundary conditions on the top/bottom boundaries, we instead extend the domain to $[0, 1] \ \times \ [-0.5, 1.5]$ and apply periodic boundary conditions on the top/bottom boundaries. As these boundaries on the extended domain are outside of the domain of influence of the shock cloud interaction over the time range of the simulation, the effect of this modified setup on the flow field is negligible but it helps alleviate any issues arising from numerical errors compounding at free boundaries. The remaining left and right boundary conditions were set as Neumann and Dirichlet, respectively, while periodicity was enforced along the $z$ direction. 

    \begin{figure}[htbp!]
        \centering
        \subfloat[Density]{\label{fig:cloud_d} \adjustbox{width=0.33\linewidth,valign=b}{\includegraphics[width=\textwidth]{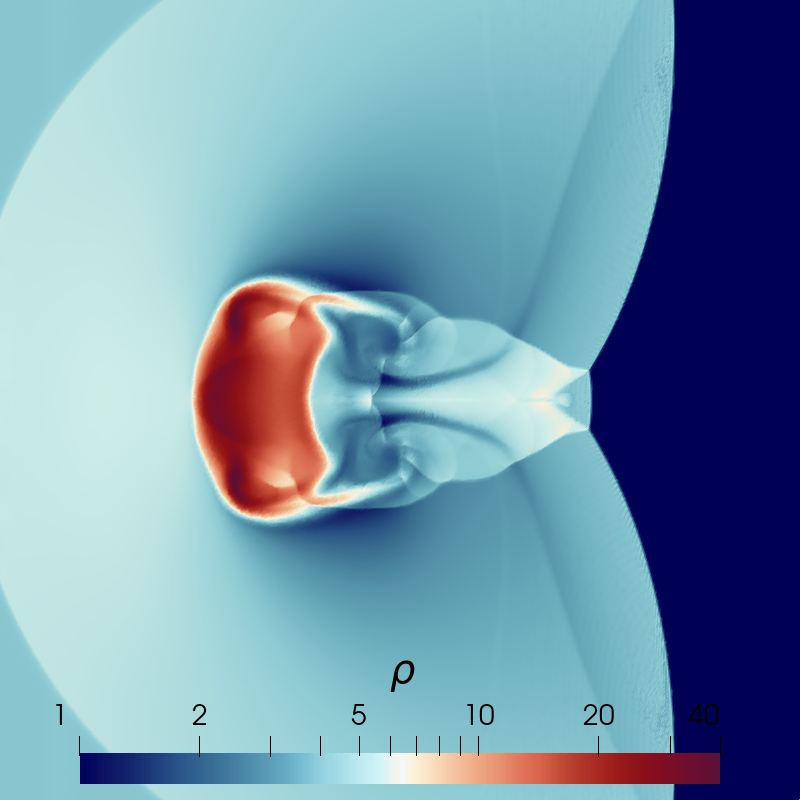}}}
        ~
        \subfloat[Pressure]{\label{fig:cloud_p} \adjustbox{width=0.33\linewidth,valign=b}{\includegraphics[width=\textwidth]{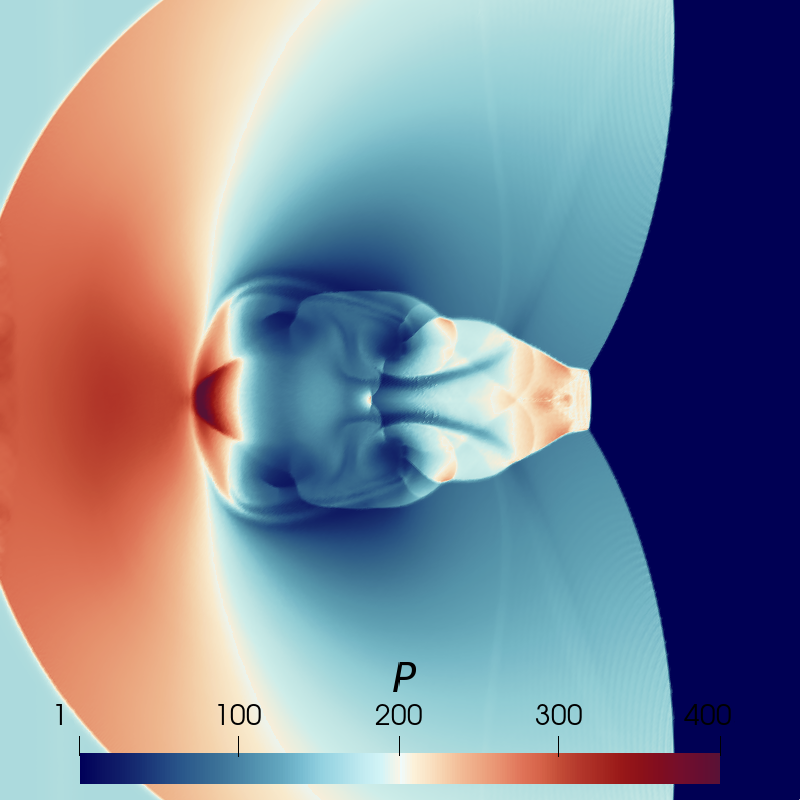}}}
        ~
        \subfloat[Magnetic pressure]{\label{fig:cloud_m} \adjustbox{width=0.33\linewidth,valign=b}{\includegraphics[width=\textwidth]{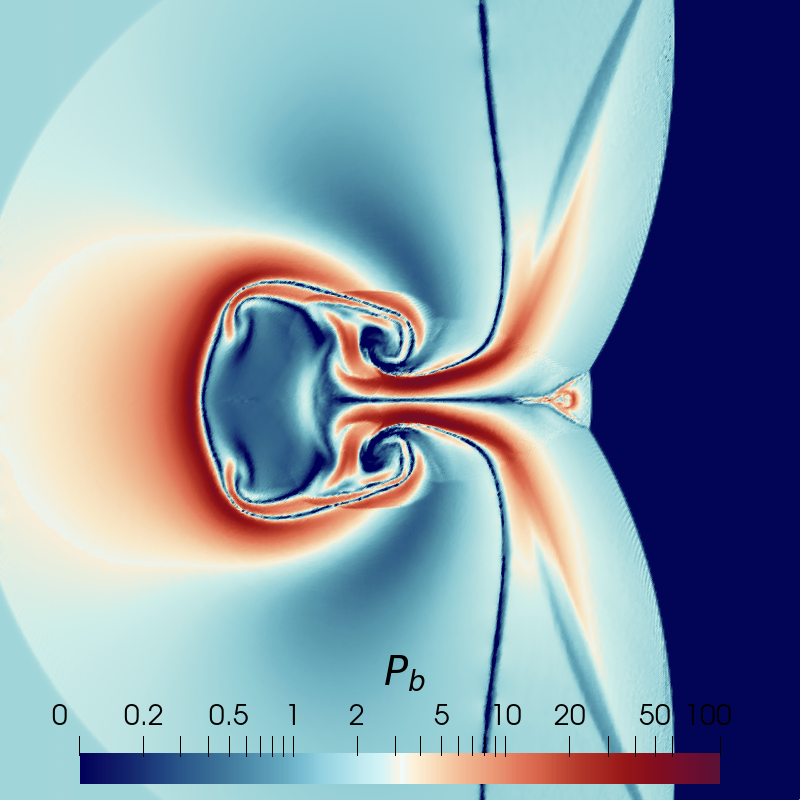}}}
        \newline
        \caption{\label{fig:cloud} Contours of density (left), pressure (middle), and magnetic pressure (right) on the subregion $[0, 1]^2$ for the shock cloud interaction problem at $t = 0.06$ computed using a $\mathbb P_2$ FR scheme with $400^2$ elements.}
    \end{figure}

\begin{figure}[htbp!]
        \centering
        \subfloat[Divergence]{\adjustbox{width=0.33\linewidth,valign=b}{\includegraphics[width=\textwidth]{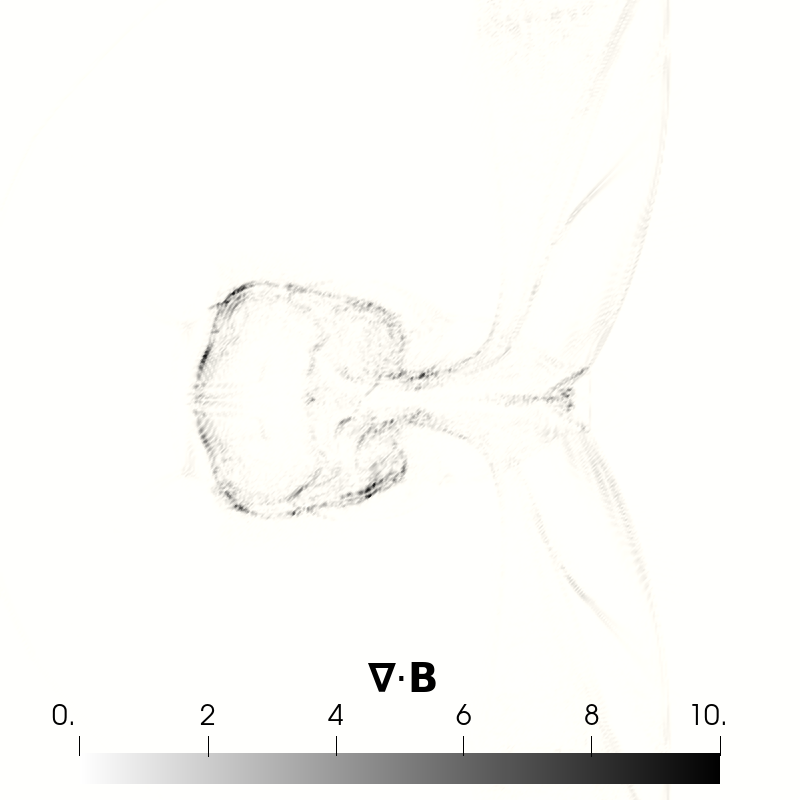}}}
        ~
        \subfloat[Limiting factor]{\adjustbox{width=0.33\linewidth,valign=b}{\includegraphics[width=\textwidth]{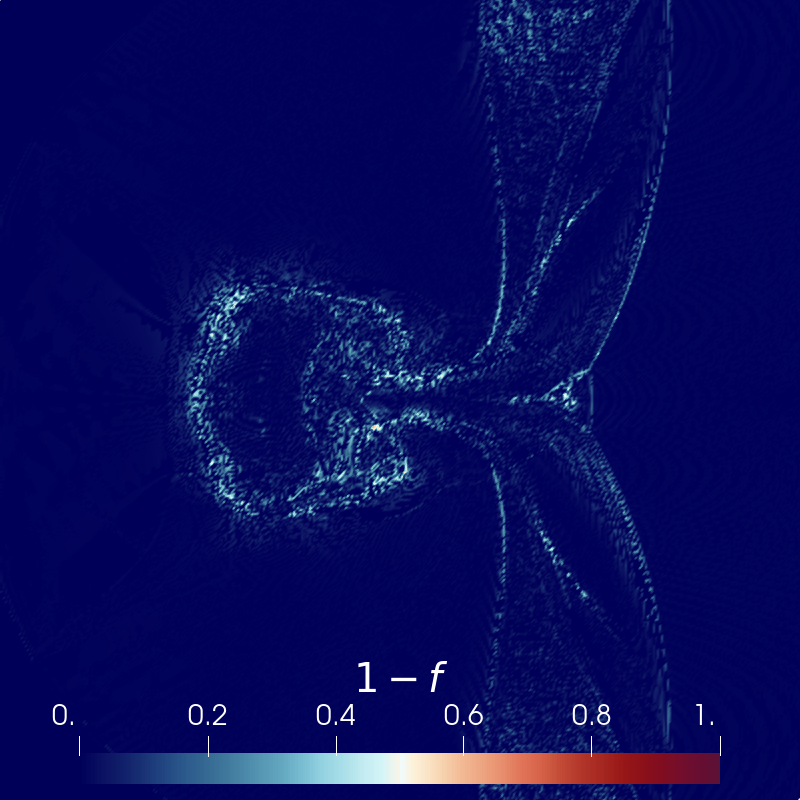}}}
        \caption{\label{fig:cloud_divbf} 
        Contours of absolute magnetic divergence (left) and limiting factor (right) on the subregion $[0, 1]^2$ for the shock cloud interaction problem at $t = 0.06$ computed using a $\mathbb P_2$ FR scheme with $400^2$ elements.
        }
    \end{figure}

To perform a comparison of the proposed approach to a third-order DG scheme augmented with a WENO limiter presented in \citet{Wu2018}, an identical problem setup is used with a $\mathbb P_2$ scheme on $400^2$ mesh (with respect to the original domain size $\Omega = [0,1]^2$). The predicted contours of density, pressure, and magnetic pressure at $t = 0.06$ computed using a time step of $\Delta t = 4{\cdot}10^{-6}$ are shown in \cref{fig:cloud}. The results show good resolution of the strong discontinuities in the various fields without any observable spurious oscillations. Furthermore, small-scale features in the cloud region of the density and magnetic fields were not excessively dissipated, and the symmetry of the problem was well-preserved. Additionally, the contours of the absolute magnetic divergence and the limiting factor are presented in \cref{fig:cloud_divbf}, which indicate that regions of strong flow gradients (i.e., discontinuities) tend to show the largest divergence errors as well as the most limiting. 

\begin{figure}[htbp!]
        \centering
        \subfloat[Entropy filter]{\adjustbox{width=0.33\linewidth,valign=b}{\includegraphics[width=\textwidth]{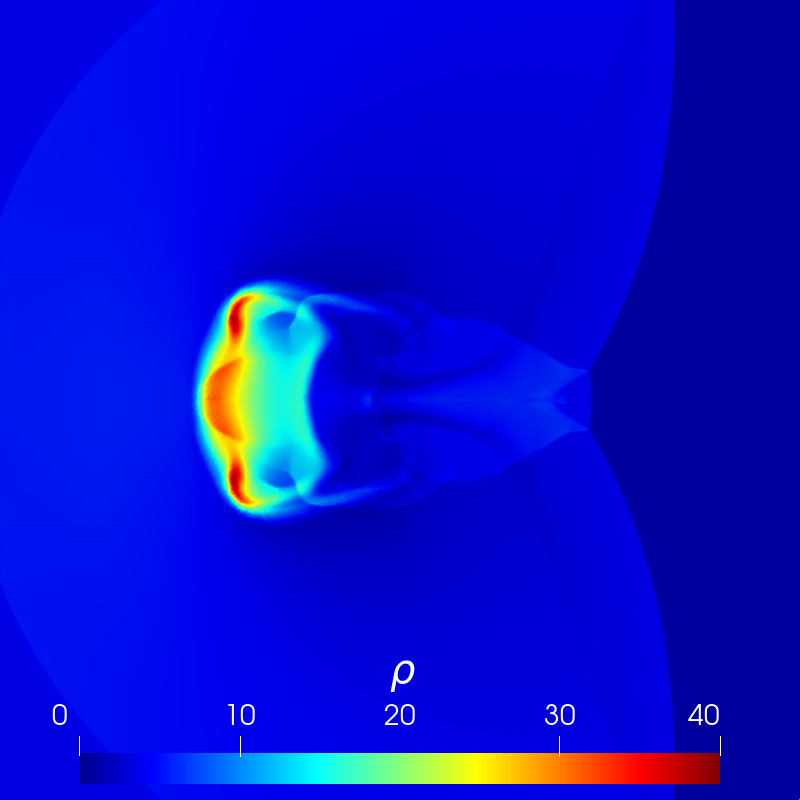}}}
        ~
        \subfloat[DG with WENO limiter]{\adjustbox{width=0.33\linewidth,valign=b}{\includegraphics[width=\textwidth]{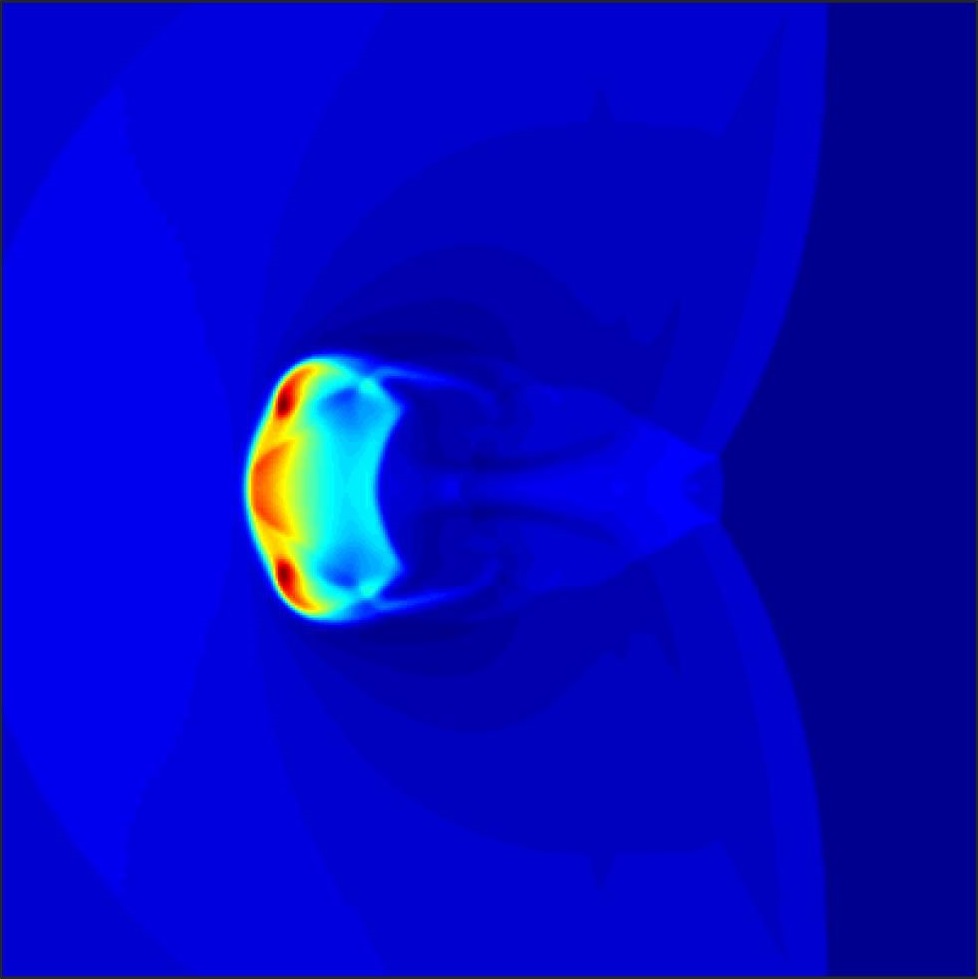}}}
        \caption{\label{fig:cloud_comp} 
        Comparison of the contours of density computed by the proposed entropy filtering approach (left) and the positivity-preserving DG scheme augmented with a WENO limiter of \citet{Wu2018} (right).
        }
    \end{figure}

A comparison to the method of \citet{Wu2018} (Fig. 1) is shown in \cref{fig:cloud_comp}. Note that some discrepancy in the color schemes between the two images may be present. The proposed scheme was roughly equally performative in terms of the resolution of discontinuities and marginally better at resolving small-scale flow features on the trailing side of the cloud. Without the positivity-preserving filtering approach, the scheme diverged due to negative pressure in the solution.

\subsection{Magnetized blast}
    \begin{figure}[htbp!]
        \centering
        \subfloat[Density]{
        \adjustbox{width=0.4\linewidth,valign=b}{\includegraphics[width=\textwidth]{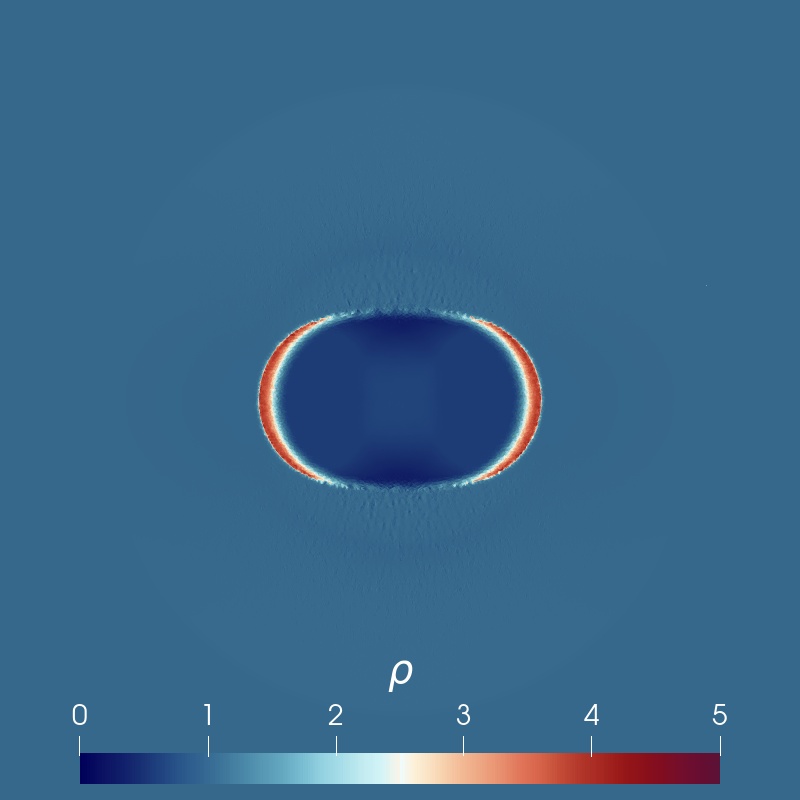}}}
        ~
        \subfloat[Velocity Magnitude]{
        \adjustbox{width=0.4\linewidth,valign=b}{\includegraphics[width=\textwidth]{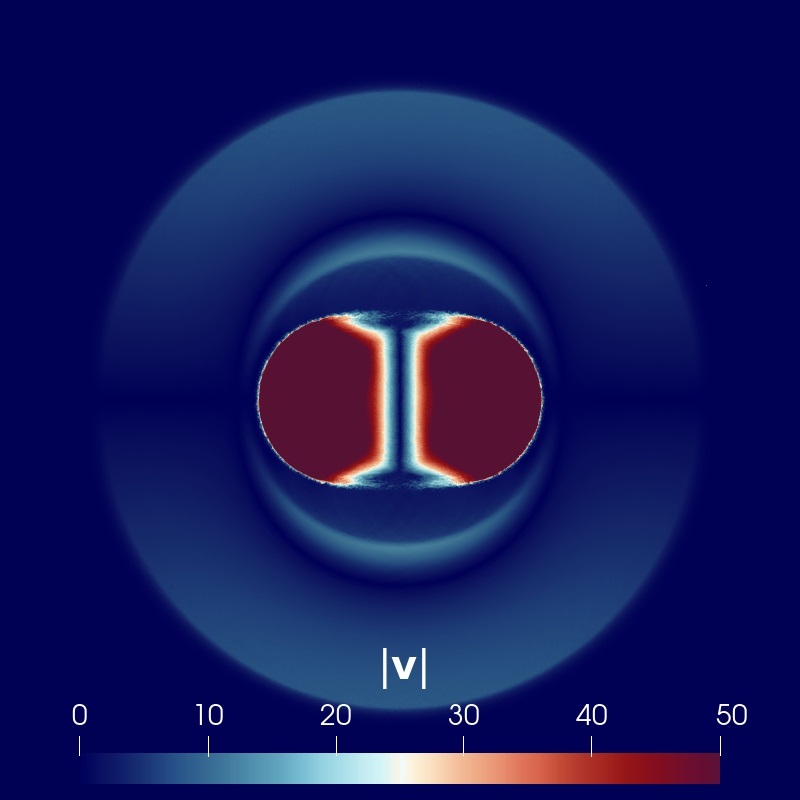}}}
        \newline
        \subfloat[Pressure]{
        \adjustbox{width=0.4\linewidth,valign=b}{\includegraphics[width=\textwidth]{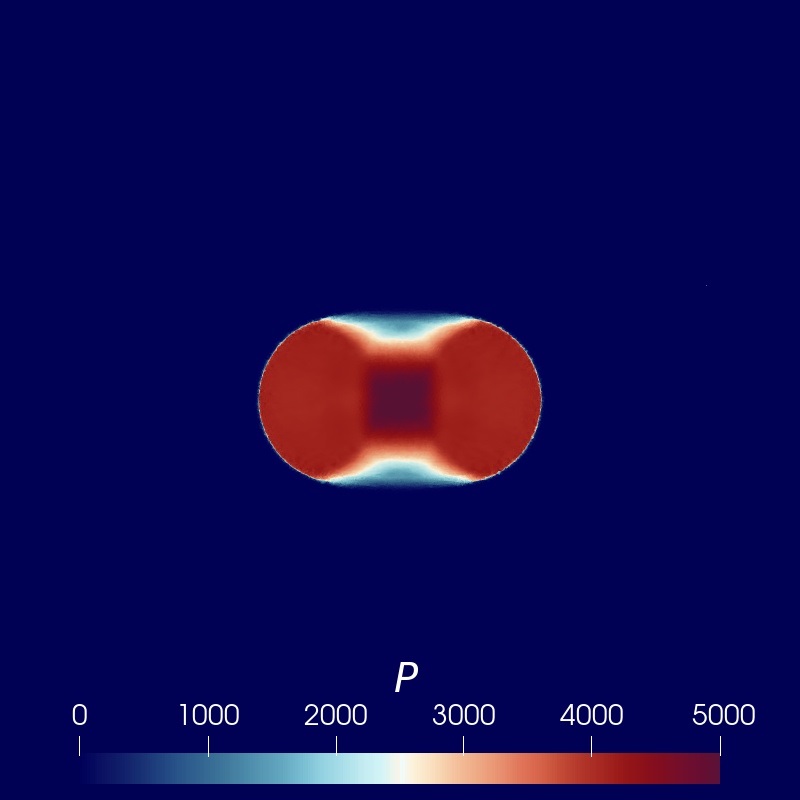}}}
        ~
        \subfloat[Magnetic Pressure]{
        \adjustbox{width=0.4\linewidth,valign=b}{\includegraphics[width=\textwidth]{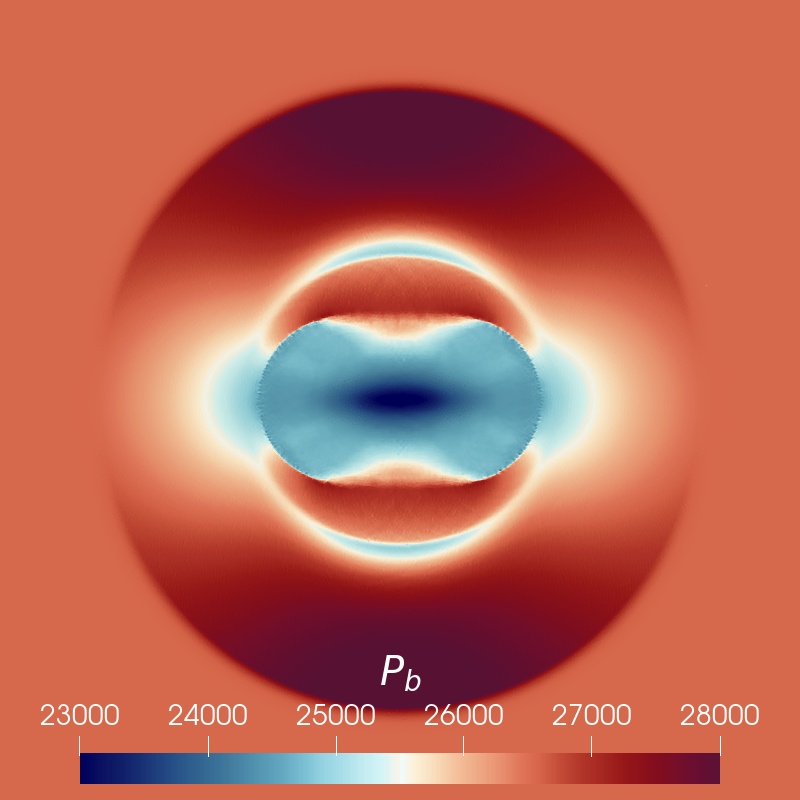}}}
        \newline
        \caption{\label{fig:blast_coarse} Contours of density (top left), velocity magnitude (top right), pressure (bottom left), and magnetic pressure (bottom right) for the magnetized blast problem at $t = 0.001$ computed using a $\mathbb P_4$ scheme on an unstructured mesh with an average edge length of $h = 1/200$.}
    \end{figure}
As a stress test for the positivity-preserving property of the proposed scheme for extreme flow conditions, a modified form of the magnetized blast wave problem of \citet{Zachary1994} and \citet{Balsara1999} was considered. In this problem, a blast wave is driven by a spherical overpressure region in the center of the domain surrounded by a low plasma-beta ambient state, resulting in strong magnetosonic shocks. The problem is solved on the periodic domain $\Omega = [-0.5, 0.5]^2$, and the initial conditions are given as
\begin{equation}
    \mathbf{q}(\mathbf{x}, 0) =    
    \begin{bmatrix} 
       \rho\\
        u\\
        v\\
        B_x\\
        B_y\\
        P
    \end{bmatrix}
    = \begin{cases}
    \mathbf{q}_e, \quad \mathrm{if} \ \sqrt{x^2 + y^2} \leq 0.1, \\
    \mathbf{q}_a, \quad \mathrm{else},
    \end{cases}
    \quad \mathrm{where} \quad 
    \mathbf{q}_e
    =
    \begin{bmatrix}
        1 \\
        0 \\
        0 \\
        B_0 \\
        0 \\
        P_e
    \end{bmatrix}
    \quad \mathrm{and} \quad \mathbf{q}_a
    =
    \begin{bmatrix}
        1 \\
        0 \\
        0 \\
        B_0 \\
        0 \\
        P_a
    \end{bmatrix}.
\end{equation}
The specific heat ratio is set as $5/3$. While the original problem setup uses $P_a = 0.1$, $P_e = 10^3$, and $B_0 = 100/\sqrt{4\pi}$, we consider the much more extreme case presented in \citet{Wu2018} with $P_a = 0.1$, $P_e = 10^4$, and $B_0 = 1000/\sqrt{4\pi}$, resulting in a very large pressure ratio of $10^5$ and a very small plasma-beta of $\beta \approx 2.5{\cdot}10^{-4}$. As these conditions are quite extreme, the scheme would diverge almost instantly in the absence of any positivity-preserving modifications. 

    \begin{figure}[htbp!]
        \centering
        \subfloat[Density]{
        \adjustbox{width=0.4\linewidth,valign=b}{\includegraphics[width=\textwidth]{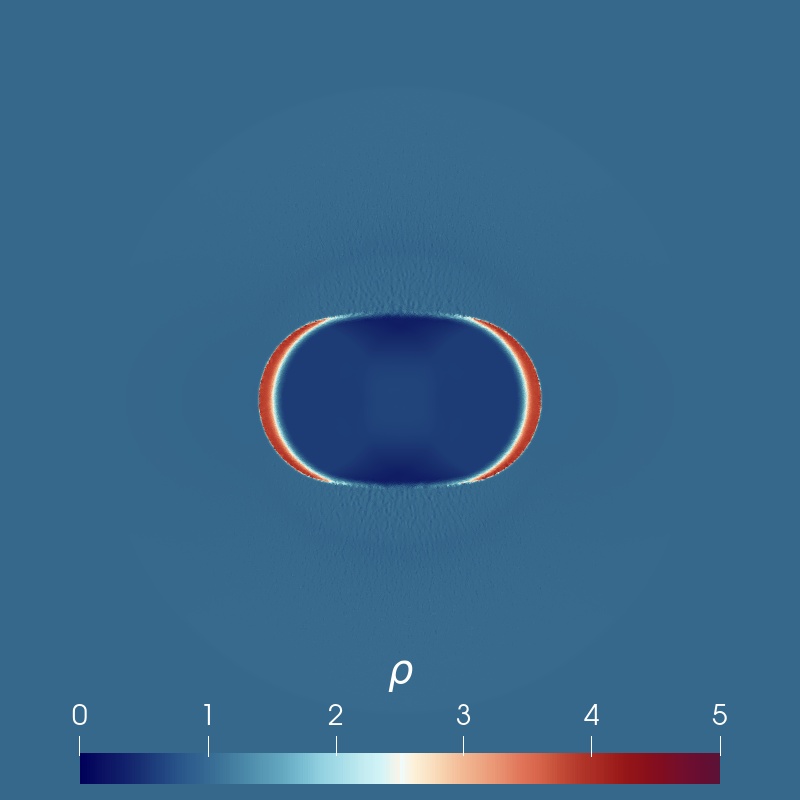}}}
        ~
        \subfloat[Velocity Magnitude]{
        \adjustbox{width=0.4\linewidth,valign=b}{\includegraphics[width=\textwidth]{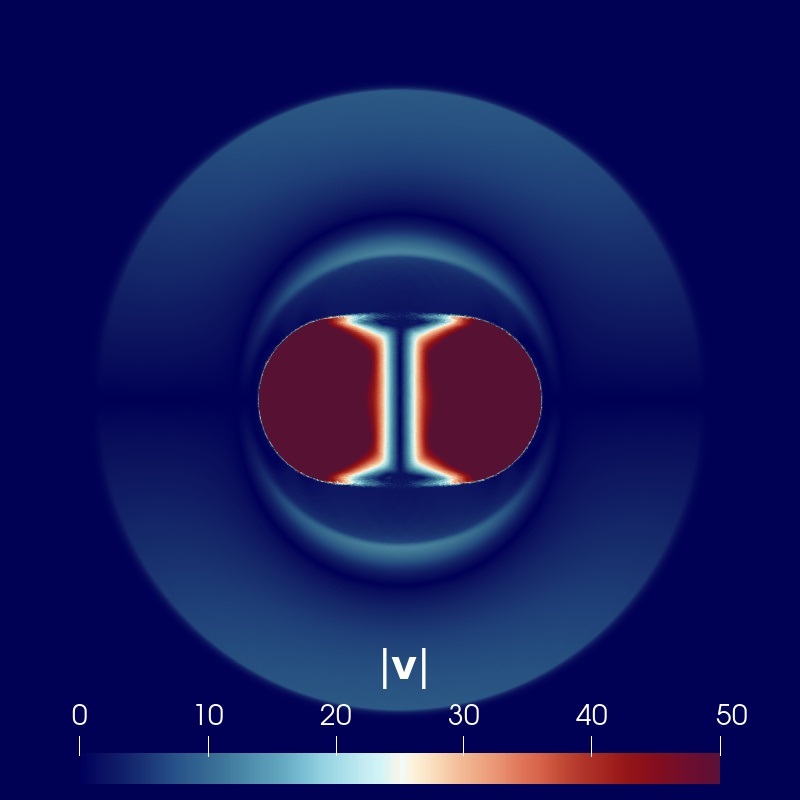}}}
        \newline
        \subfloat[Pressure]{
        \adjustbox{width=0.4\linewidth,valign=b}{\includegraphics[width=\textwidth]{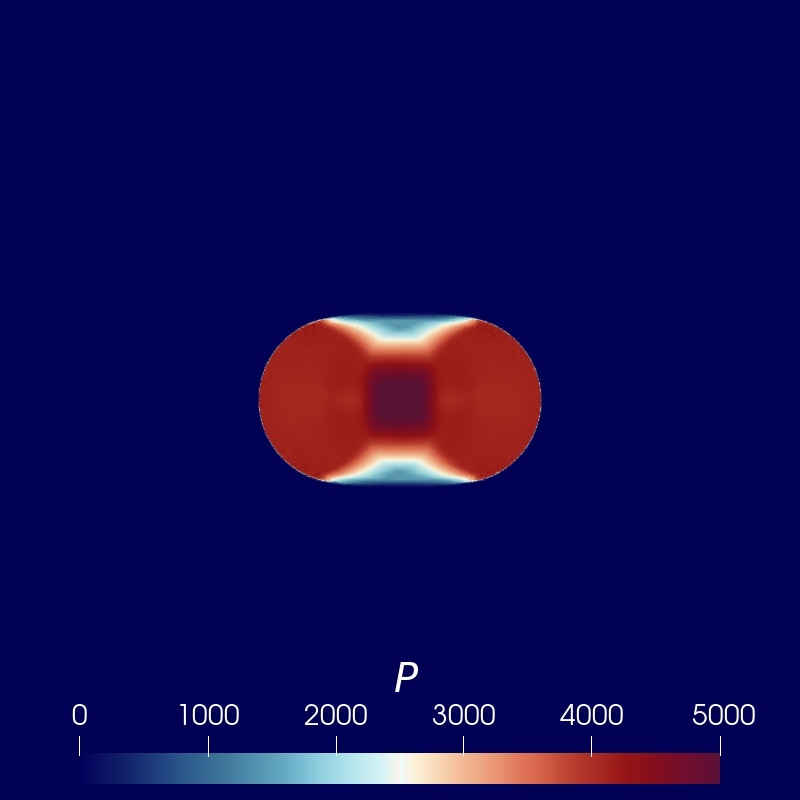}}}
        ~
        \subfloat[Magnetic Pressure]{
        \adjustbox{width=0.4\linewidth,valign=b}{\includegraphics[width=\textwidth]{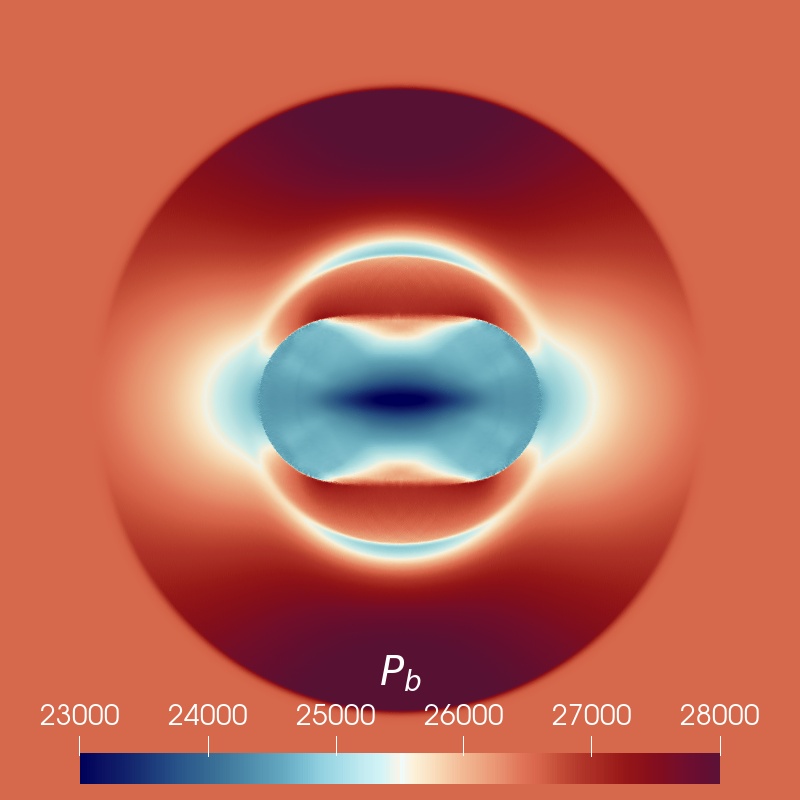}}}
        \newline
        \caption{\label{fig:blast_fine} Contours of density (top left), velocity magnitude (top right), pressure (bottom left), and magnetic pressure (bottom right) for the magnetized blast problem at $t = 0.001$ computed using a $\mathbb P_4$ scheme on an unstructured mesh with an average edge length of $h = 1/400$.}
    \end{figure}

\begin{figure}[htbp!]
        \centering
        \subfloat[Coarse mesh]{
        \adjustbox{width=0.4\linewidth,valign=b}{\includegraphics[width=\textwidth]{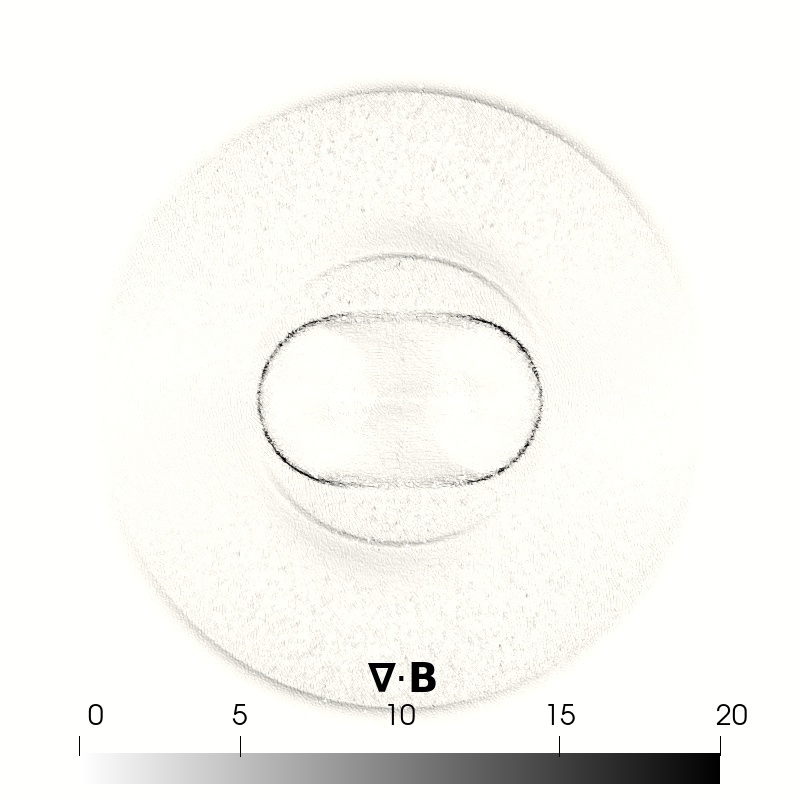}}}
        ~
        \subfloat[Fine mesh]{
        \adjustbox{width=0.4\linewidth,valign=b}{\includegraphics[width=\textwidth]{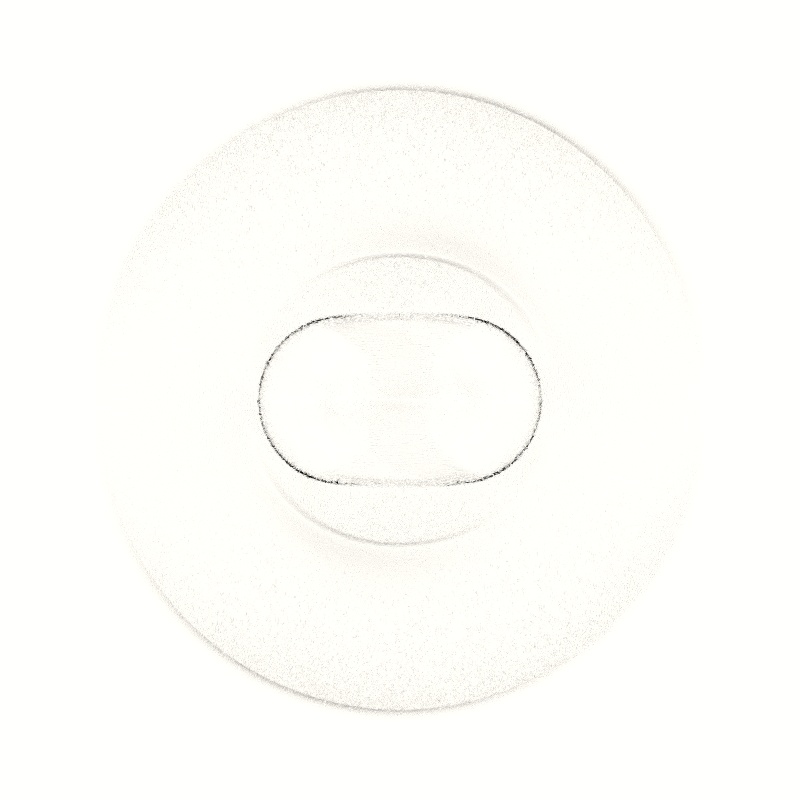}}}
        \newline
        \caption{\label{fig:blast_divb} Contours of absolute magnetic divergence for the magnetized blast problem at $t = 0.001$ computed using a $\mathbb P_4$ scheme on the coarse mesh (left) and fine mesh (right).}
    \end{figure}

To verify that the proposed scheme can be easily extended to unstructured grids, the problem was solved on triangular meshes using a $\mathbb P_4$ scheme. A coarse mesh and a fine mesh were generated, consisting of approximately $1.2{\cdot}10^5$ elements with an average edge length of $h = 1/200$ and approximately $5{\cdot}10^5$ elements with an average edge length of $h = 1/400$, respectively. For this case, the HLL Riemann solver was used as it was found to be much better behaved in these extreme conditions than the HLLC Riemann solver, although both approaches were properly stabilized with the proposed entropy filtering method. The contours of density, velocity magnitude, pressure, and magnetic pressure at $t = 0.001$ are shown in \cref{fig:blast_coarse} and \cref{fig:blast_fine} for the coarse and fine meshes, respectively, computed using time steps of $\Delta t = 2{\cdot}10^{-7}$ and $\Delta t = 1{\cdot}10^{-7}$. Even with these extreme conditions on unstructured grids, the predicted solutions were well-behaved, and both the coarse and fine meshes showed excellent resolution of the various discontinuities in the velocity, pressure, and magnetic fields with sub-element resolution and no observable spurious oscillations. Furthermore, the numerical width of the discontinuities decreased appropriately with increasing resolution. For the density field, minor spurious oscillations were observed, primarily at lower resolution and somewhat indicative of mesh imprinting, but the strength and distribution of these oscillations decreased with the finer mesh. These observations are consistent with the case of the Brio--Wu shock tube where the density field was marginally less well-behaved than the other fields. However, the predicted fields were still very good given such extreme conditions and an unstructured mesh, indicating that the proposed approach remains robust and accurate for such flows. The contours of the absolute magnetic divergence are additionally shown in \cref{fig:blast_divb} for both the coarse and fine mesh. Similarly to the previous experiments, the divergence error was predominantly focused in regions of strong gradients, and this divergence error tended to decrease with increasing mesh resolution. 

\begin{figure}[htbp!]
        \centering
        \subfloat[Coarse mesh]{
        \adjustbox{width=0.4\linewidth,valign=b}{\includegraphics[width=\textwidth]{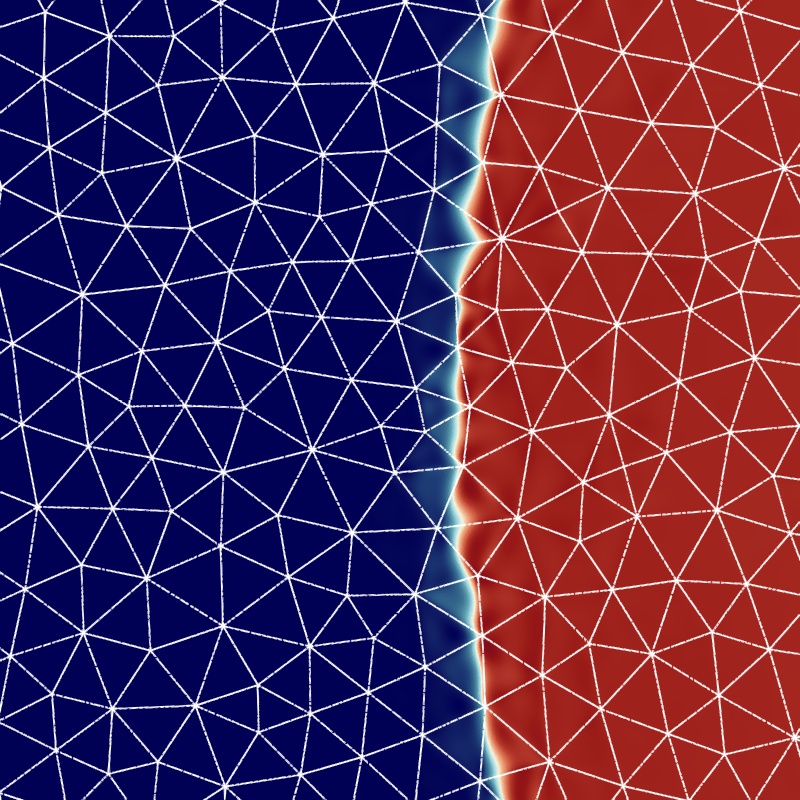}}}
        ~
        \subfloat[Fine mesh]{
        \adjustbox{width=0.4\linewidth,valign=b}{\includegraphics[width=\textwidth]{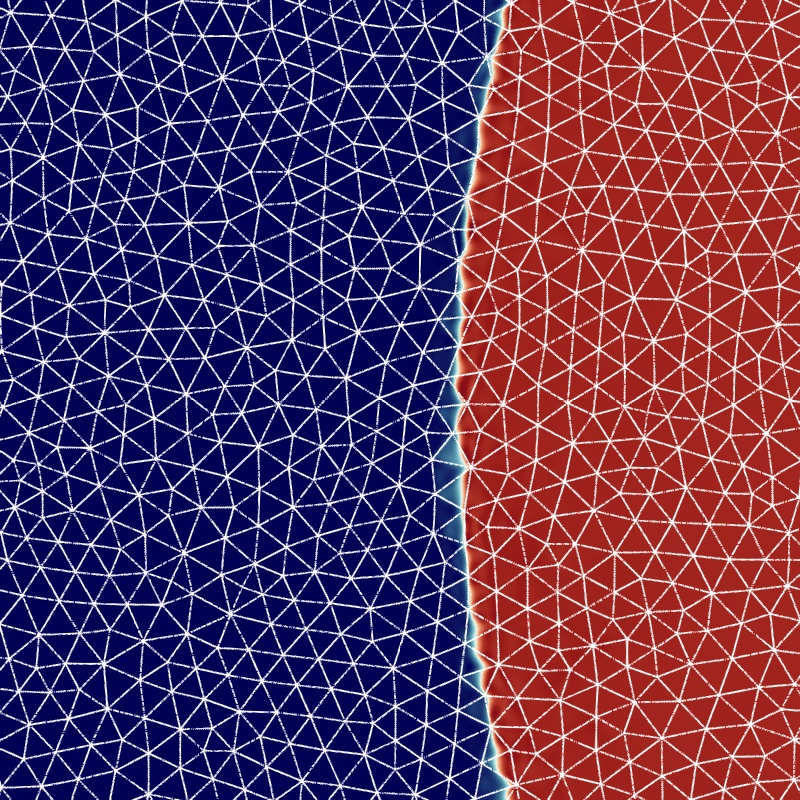}}}
        \newline
        \caption{\label{fig:blast_mesh} Enlarged view of contours of pressure with mesh overlay for the magnetized blast problem at $t = 0.001$ computed using a $\mathbb P_4$ scheme on the coarse mesh (left) and fine mesh (right). Contour scale identical to \cref{fig:blast_coarse} and \cref{fig:blast_fine}.}
    \end{figure}

To demonstrate the sub-element shock-resolving ability of the entropy filtering approach on unstructured grids, an enlarged view of the contours of pressure with the mesh overlaid is shown for the two meshes in \cref{fig:blast_mesh}. It can be seen that the shock was resolved well within the element, with the majority of the feature resolved across 1-2 solution nodes. Furthermore, this behavior persisted when the mesh resolution was increased, such that the discrete shock thickness decreased proportionally. Given the unstructured nature of the mesh and the resulting solution point distribution, the circular shape of the shock front was still well-represented, with even better approximation at higher mesh resolutions. These results indicate that the proposed approach can be extended to unstructured grids in a straightforward manner without appreciably sacrificing its efficiency or performance at resolving discontinuities. 

\subsection{Three-dimensional Rayleigh--Taylor instability}
A final evaluation of the proposed approach and the extension to three-dimensional flows was performed through the simulation of a magnetized three-dimensional Rayleigh--Taylor instability. The problem consists of a denser gas resting on top of a lighter gas under the effect of a gravitational field initially in equilibrium with the pressure gradient. Instabilities arise in the form of ``bubbles'' of the lighter gas rising and ``fingers'' of the heavier gas descending, after which nonlinear momentum transport drives the flow to a turbulent mixing state. The problem is solved on the domain $[-L/2, L/2]^2 \times [-L, L]$, where $L = 1$, and the initial conditions are given as
\begin{equation}
    \mathbf{q}(\mathbf{x}, 0) =     
       [\rho, 
        u, 
        v, 
        w, 
        B_x, 
        B_y,
        B_z,
        P]^T
    = \begin{cases}
    \mathbf{q}_l, \quad \, \mathrm{if} \ z \leq 0, \\
    \mathbf{q}_h, \quad \mathrm{else},
    \end{cases}
\end{equation}
where
\begin{equation}
    \mathbf{q}_l
    =
    \begin{bmatrix}
        \rho_l \\
        0 \\
        0 \\
        W(x,y,z) \\
        B_0 \\
        0 \\
        0 \\
        P(z)
    \end{bmatrix}
    \quad \mathrm{and} \quad \mathbf{q}_h
    =
    \begin{bmatrix}
        \rho_h \\
        0 \\
        0 \\
        W(x,y,z) \\
        B_0 \\
        0 \\
        0 \\
        P(z)
    \end{bmatrix},
\end{equation}
for some vertical velocity perturbation $W(x,y,z)$ and initial pressure distribution $P(z)$. Periodic boundary conditions are enforced along the transverse ($x$, $y$) directions while reflecting boundary conditions are enforced along the top and bottom boundaries. A gravitational field is added to the problem, given in the form of a source term as 
\begin{equation}
    \mathbf{S_G} = -[0, 0, 0, \rho g, 0, 0, 0, \rho w g],
\end{equation}
where $g = 1$. The treatment of this gravitational field in the context of the entropy filter is taken simply as an additional term in the source term of Powell's method, and as it not stiff for this problem, it does not appreciably affect the time step restrictions of the scheme. 

The parameters for the problem are taken similarly to a scaled form of the setup in \citet{Stone2007}. The densities of the light and heavy gases are taken as $\rho_l = 1$ and $\rho_h = 3$, respectively, yielding an Atwood number of $1/2$. To enforce equilibrium in the flow, the initial pressure field is taken as
\begin{equation}
    P(z) = P_0 - \rho g z,
\end{equation}
where $P_0 = 10/\gamma$ for a specific heat ratio $\gamma = 5/3$. To seed instabilities in the flow, perturbations were added in the form of a vertical velocity field component as
\begin{equation}
    W(x,y,z) = A \cos \left(\frac{\pi z}{2 L} \right) \sin \left(\frac{4 \pi x}{L}\right) \sin \left(\frac{4 \pi y}{L}\right),
\end{equation}
where $A = 0.05$. This differs from the work of \citet{Stone2007} in that the transverse distribution of the perturbations is taken as a single deterministic mode instead of randomly-generated noise. As such, the predicted flow fields are expected to differ during the linear growth regime of the instability. 

The addition of a magnetic field significantly impacts the behavior of the Rayleigh--Taylor instability as it induces a stabilizing effect on the flow. In fact, linear stability analysis presents a cutoff magnetic field value,
\begin{equation}
    B_c = \sqrt{(\rho_h - \rho_l)g L} = \sqrt{2},
\end{equation}
above which the instability is completely damped by the magnetic field \citep{Stone2007}. We consider three variations of this flow, a hydrodynamic case where $B_0 = 0$, a weakly magnetized case where $B_0 = 0.1 B_c$, and a strongly magnetized case where $B_0 = 0.5 B_c$. The term strongly magnetized is relative in the sense that the plasma-beta is still quite high, but the magnetic field can suppress almost all potential instability modes along its orientation.

\begin{figure}[htbp!]
        \centering
        \subfloat[$t = 2$]{
        \adjustbox{width=0.24\linewidth,valign=b}{\includegraphics[width=\textwidth]{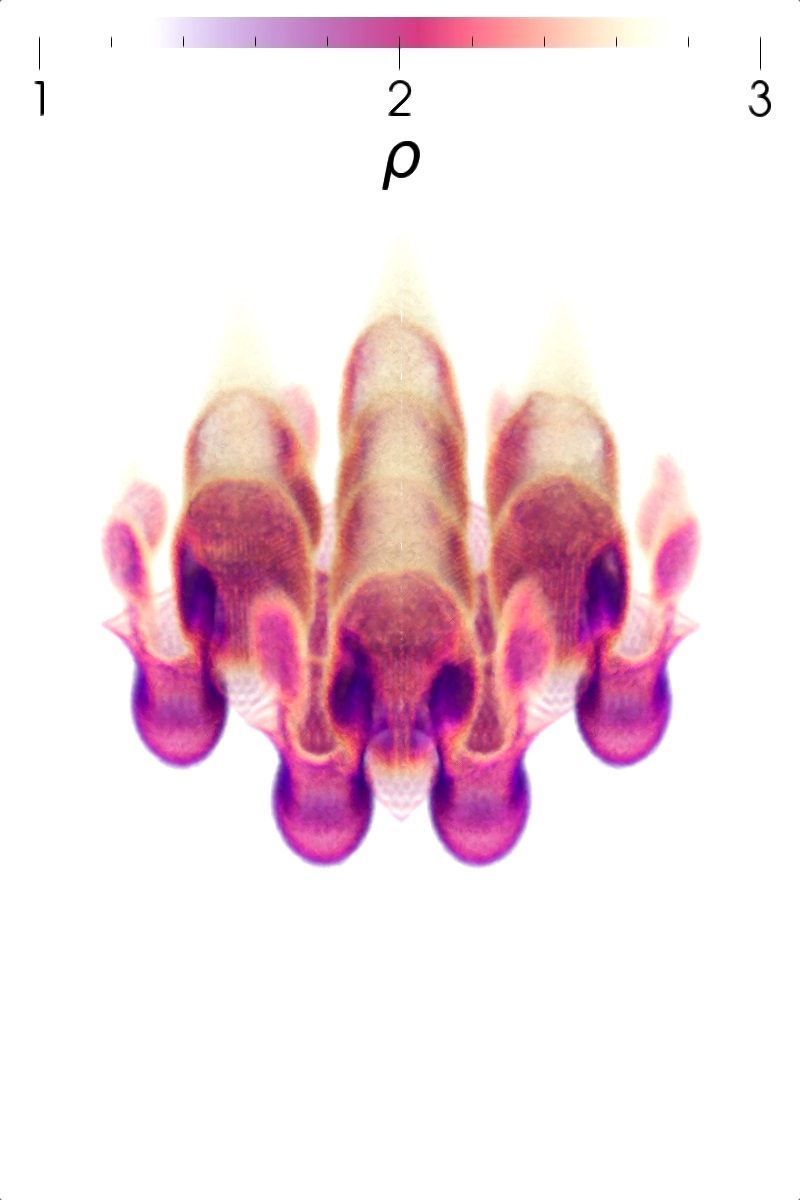}}}
        ~
        \subfloat[$t = 3$]{
        \adjustbox{width=0.24\linewidth,valign=b}{\includegraphics[width=\textwidth]{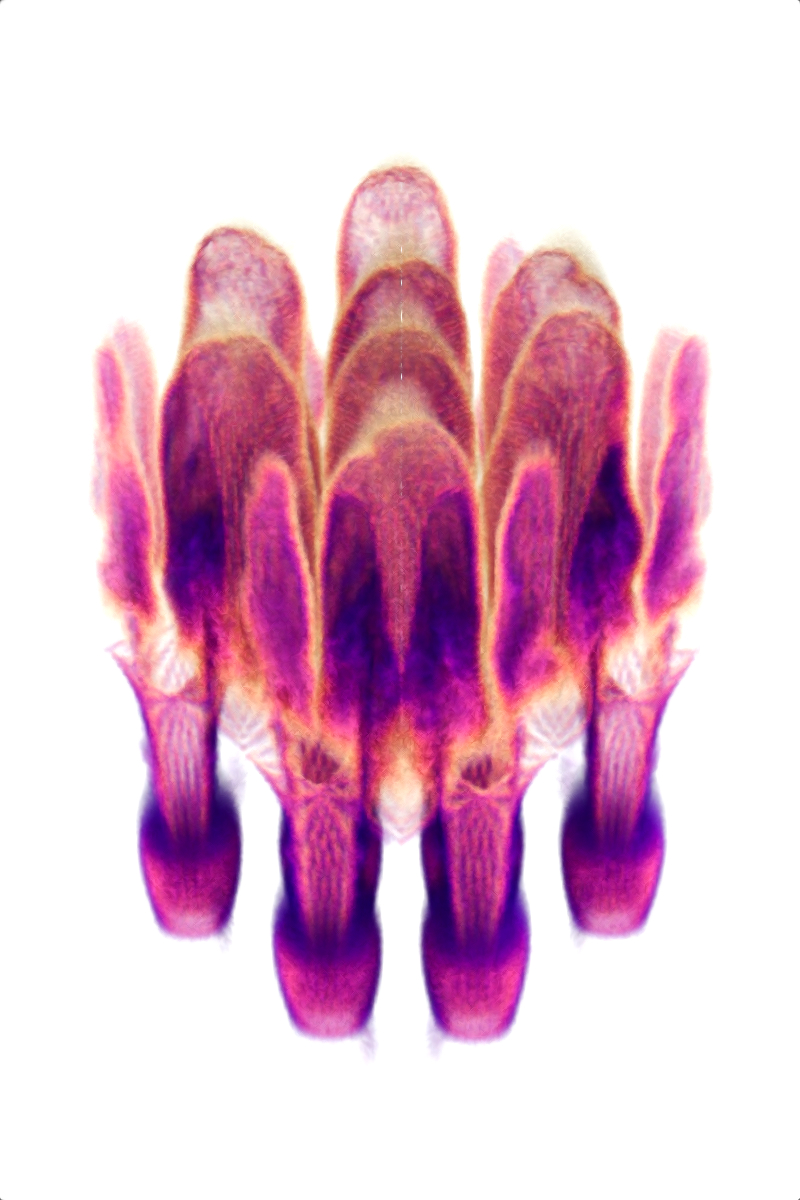}}}
        ~
        \subfloat[$t = 4$]{
        \adjustbox{width=0.24\linewidth,valign=b}{\includegraphics[width=\textwidth]{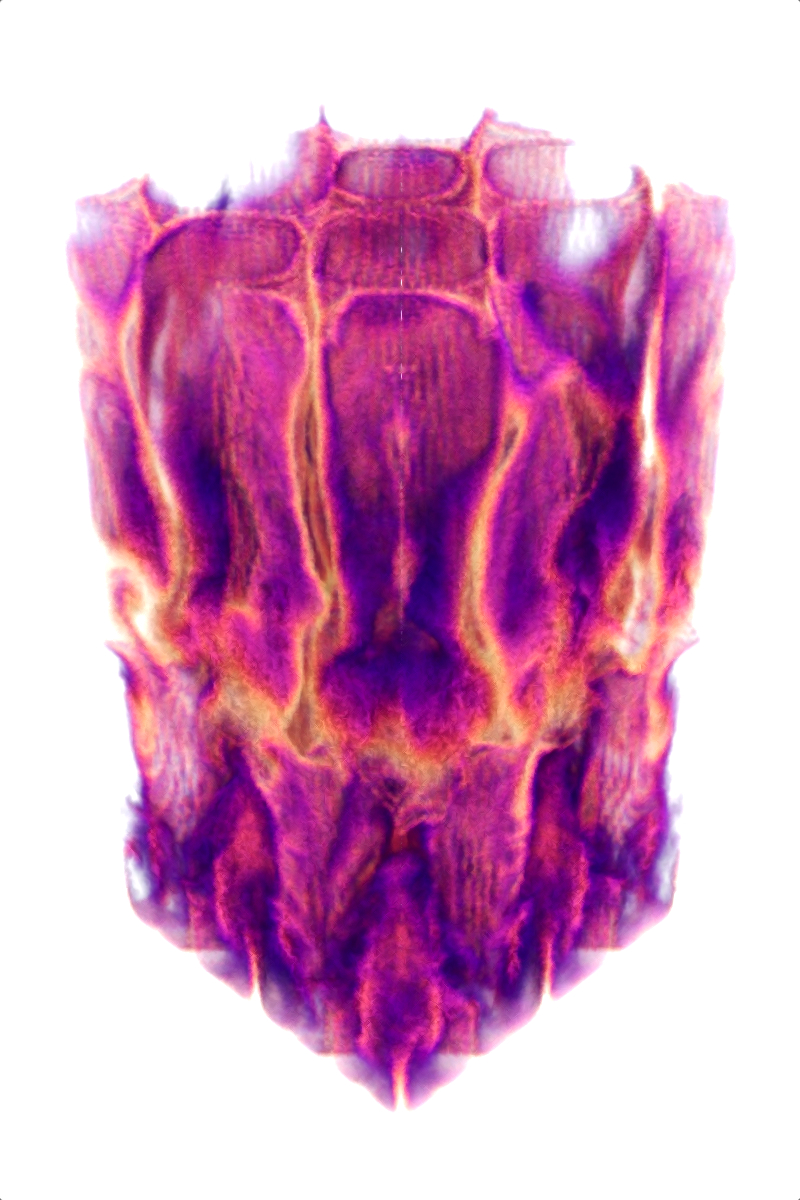}}}
        ~
        \subfloat[$t = 5$]{
        \adjustbox{width=0.24\linewidth,valign=b}{\includegraphics[width=\textwidth]{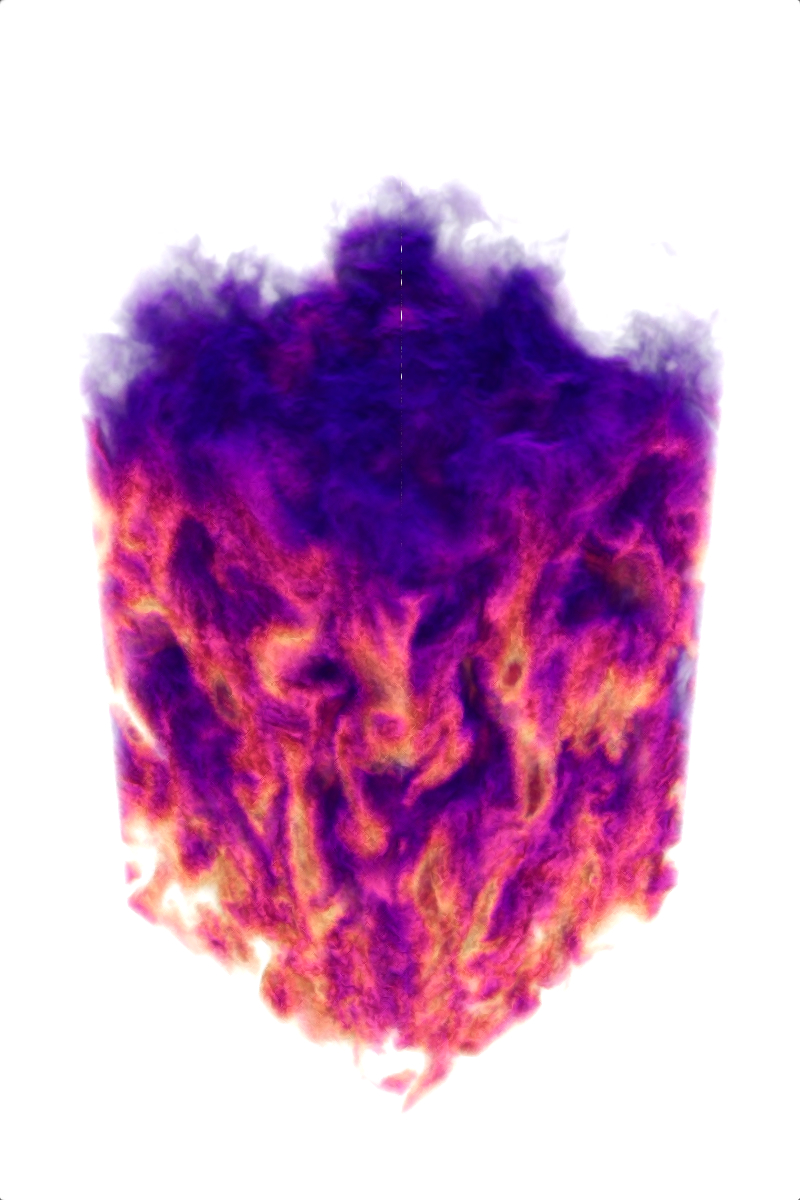}}}
        \newline
        \caption{\label{fig:rt_b0} Volume rendering of the density field for the hydrodynamic ($B_0 = 0$) Rayleigh--Taylor instability problem at varying times computed using a $\mathbb P_3$ scheme on a $N_e = 64 \times 64 \times 128$ mesh.}
    \end{figure}

\begin{figure}[htbp!]
        \centering
        \subfloat[$t = 2$]{
        \adjustbox{width=0.24\linewidth,valign=b}{\includegraphics[width=\textwidth]{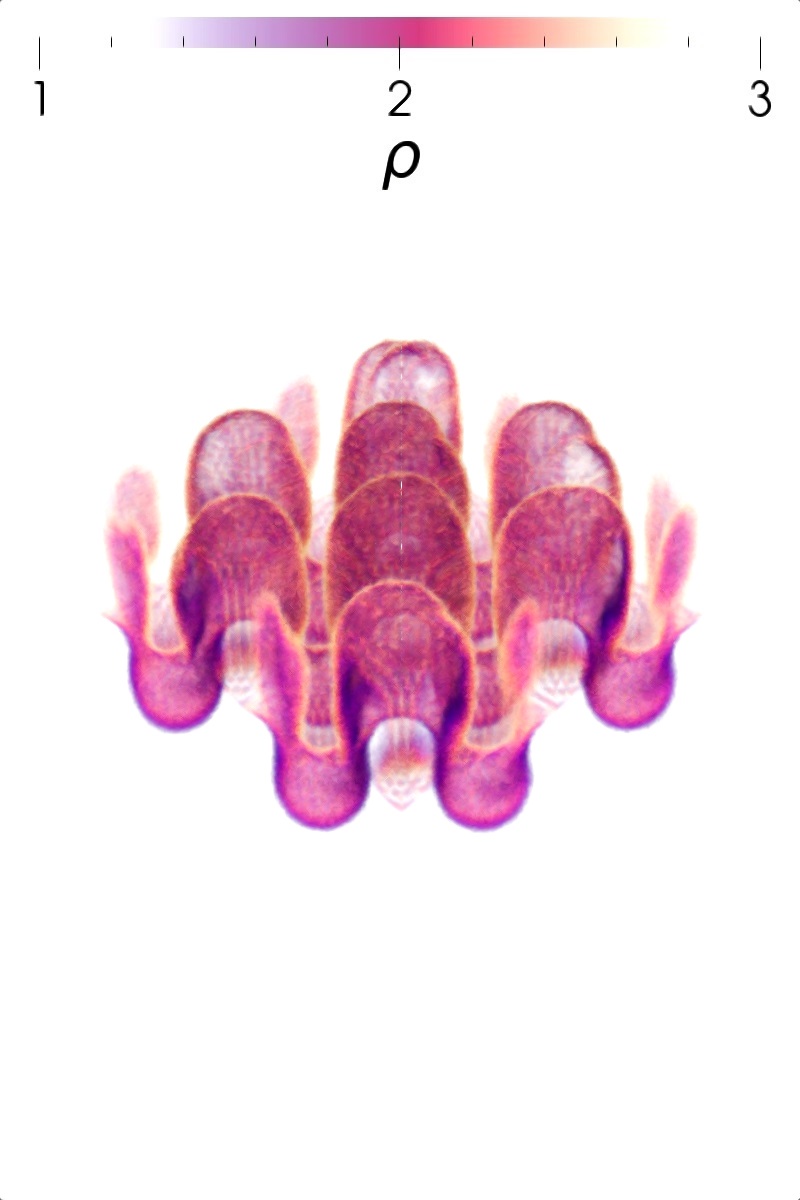}}}
        ~
        \subfloat[$t = 3$]{
        \adjustbox{width=0.24\linewidth,valign=b}{\includegraphics[width=\textwidth]{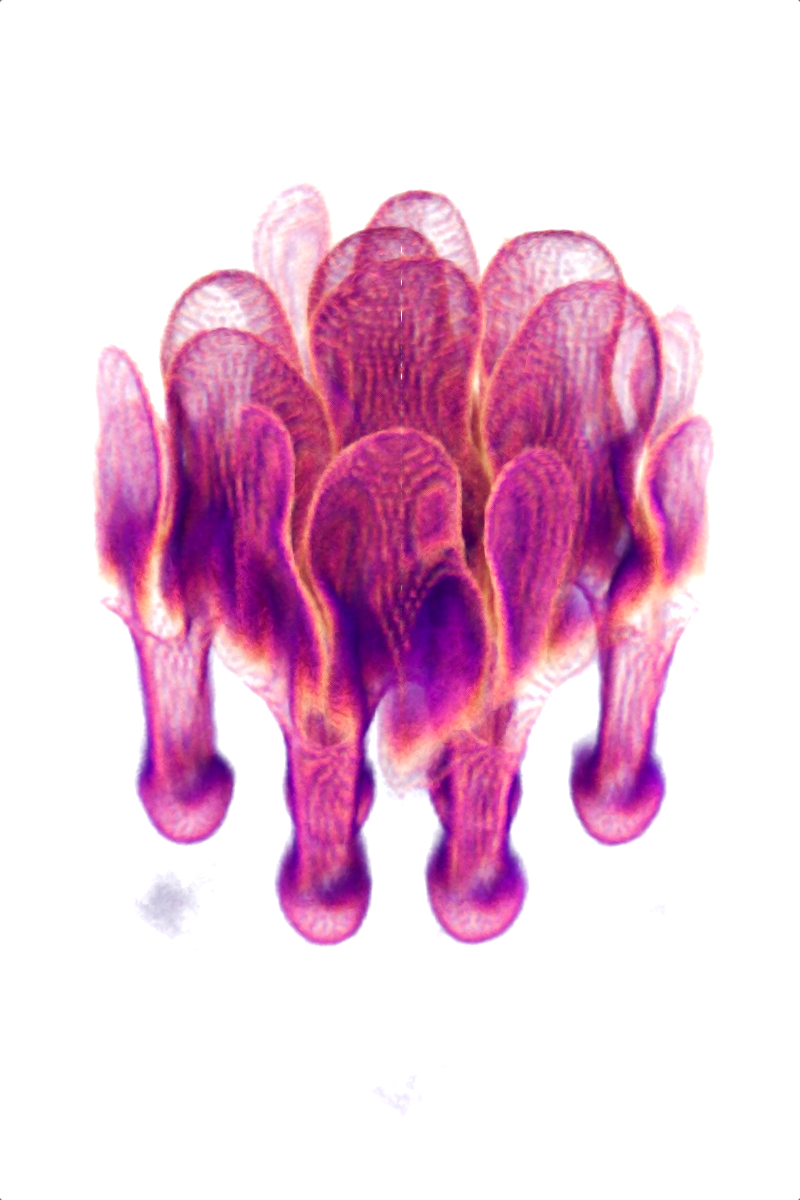}}}
        ~
        \subfloat[$t = 4$]{
        \adjustbox{width=0.24\linewidth,valign=b}{\includegraphics[width=\textwidth]{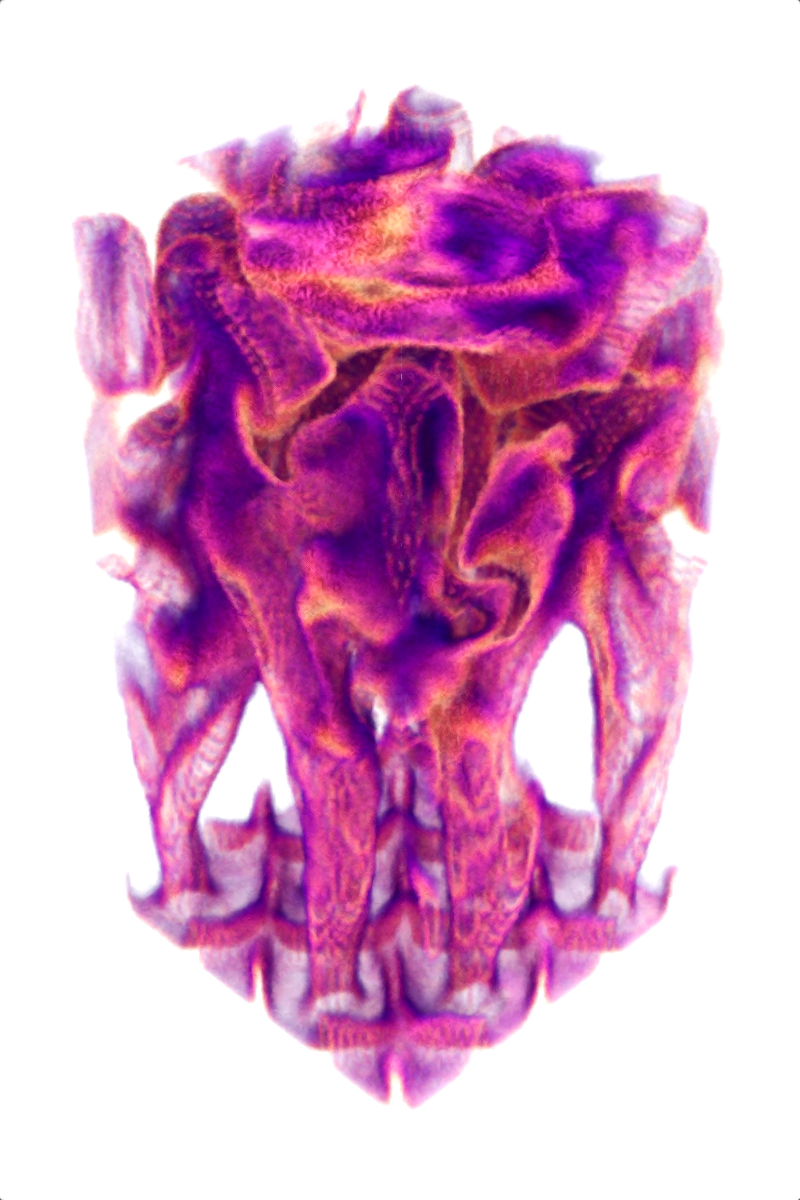}}}
        ~
        \subfloat[$t = 5$]{
        \adjustbox{width=0.24\linewidth,valign=b}{\includegraphics[width=\textwidth]{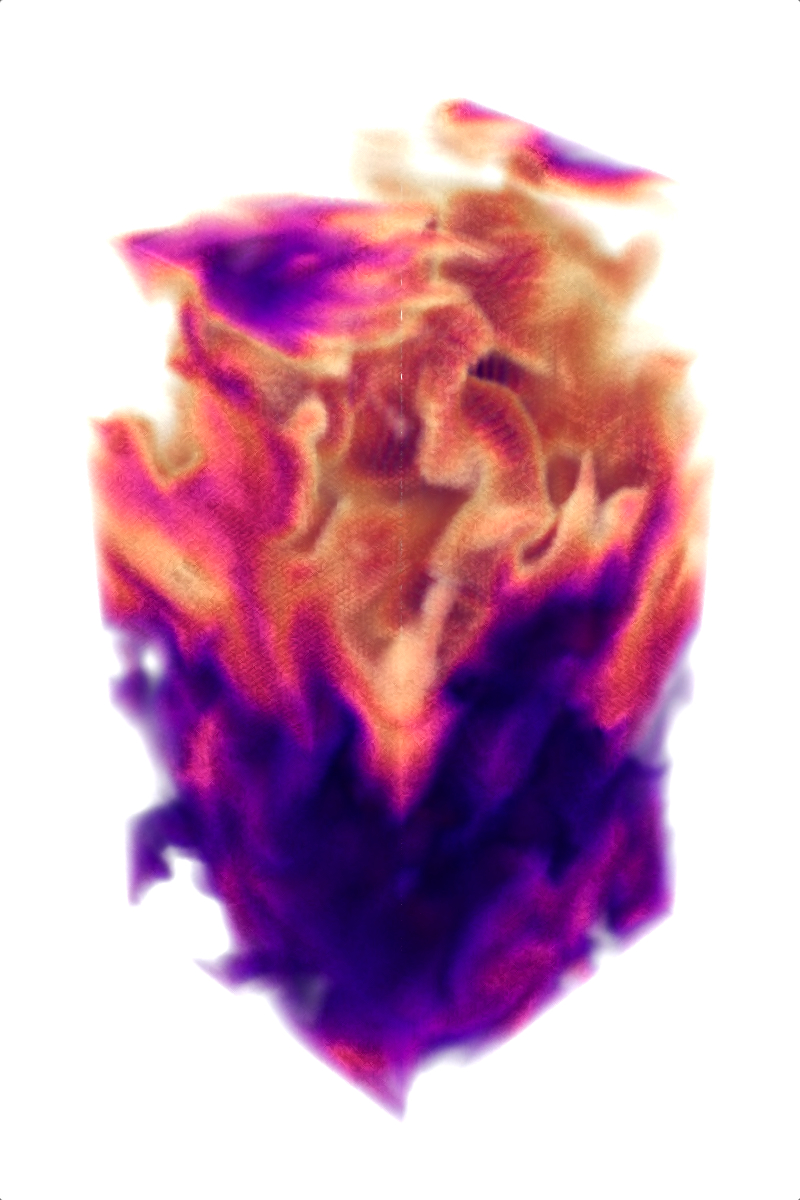}}}
        \newline
        \caption{\label{fig:rt_b14} Volume rendering of the density field for the weakly magnetized ($B_0 = 0.1 B_c$) Rayleigh--Taylor instability problem at varying times computed using a $\mathbb P_3$ scheme on a $N_e = 64 \times 64 \times 128$ mesh.}
    \end{figure}

\begin{figure}[htbp!]
        \centering
        \subfloat[$t = 2$]{
        \adjustbox{width=0.24\linewidth,valign=b}{\includegraphics[width=\textwidth]{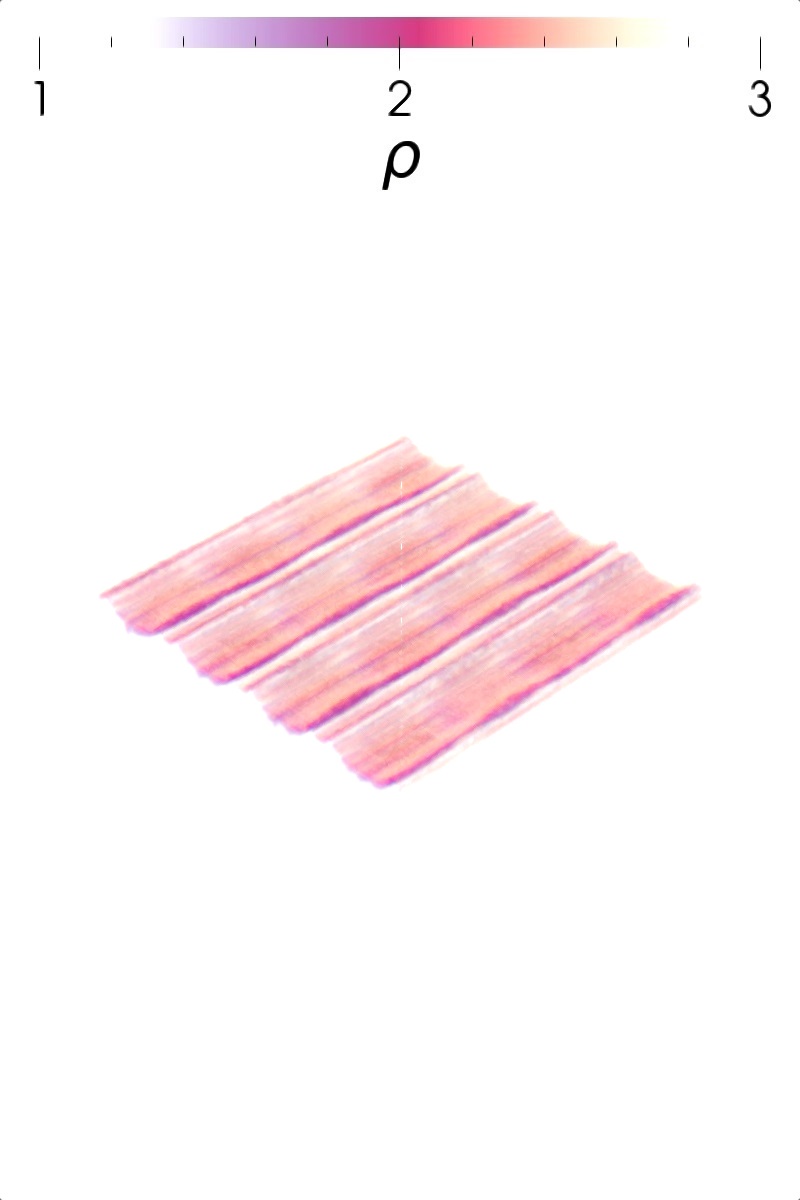}}}
        ~
        \subfloat[$t = 3$]{
        \adjustbox{width=0.24\linewidth,valign=b}{\includegraphics[width=\textwidth]{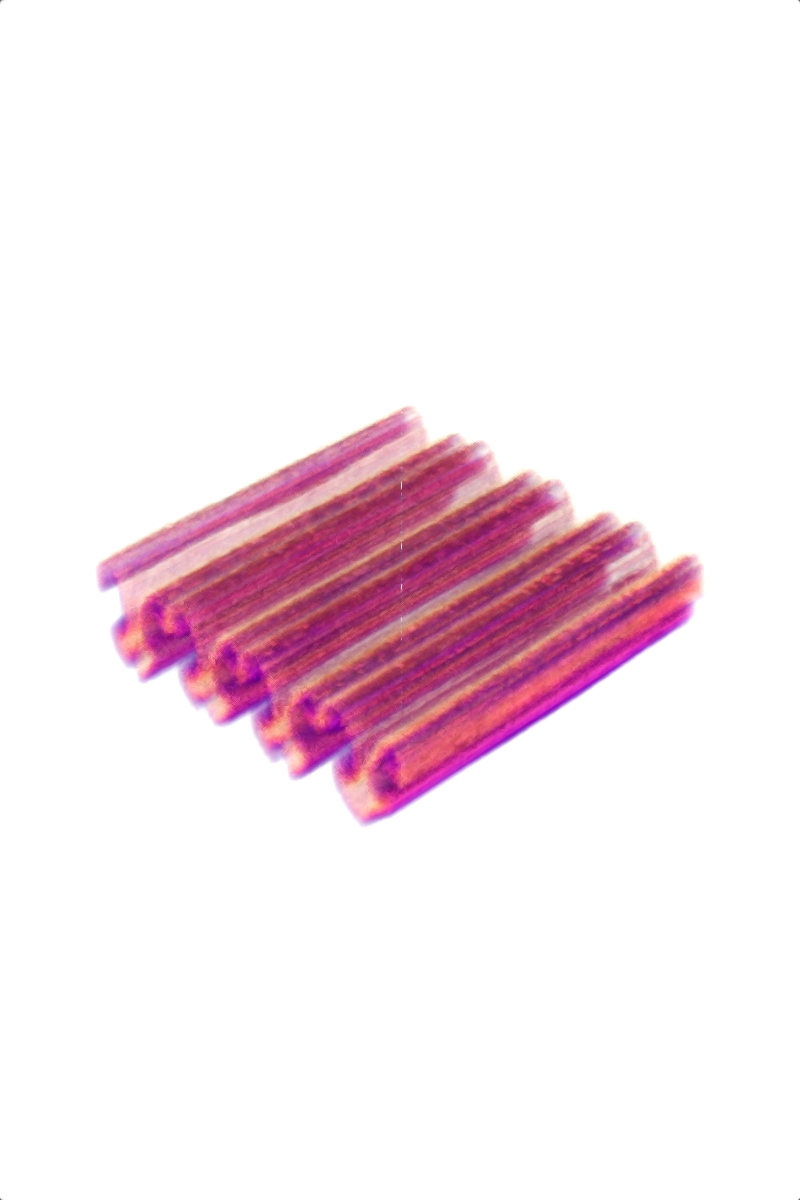}}}
        ~
        \subfloat[$t = 4$]{
        \adjustbox{width=0.24\linewidth,valign=b}{\includegraphics[width=\textwidth]{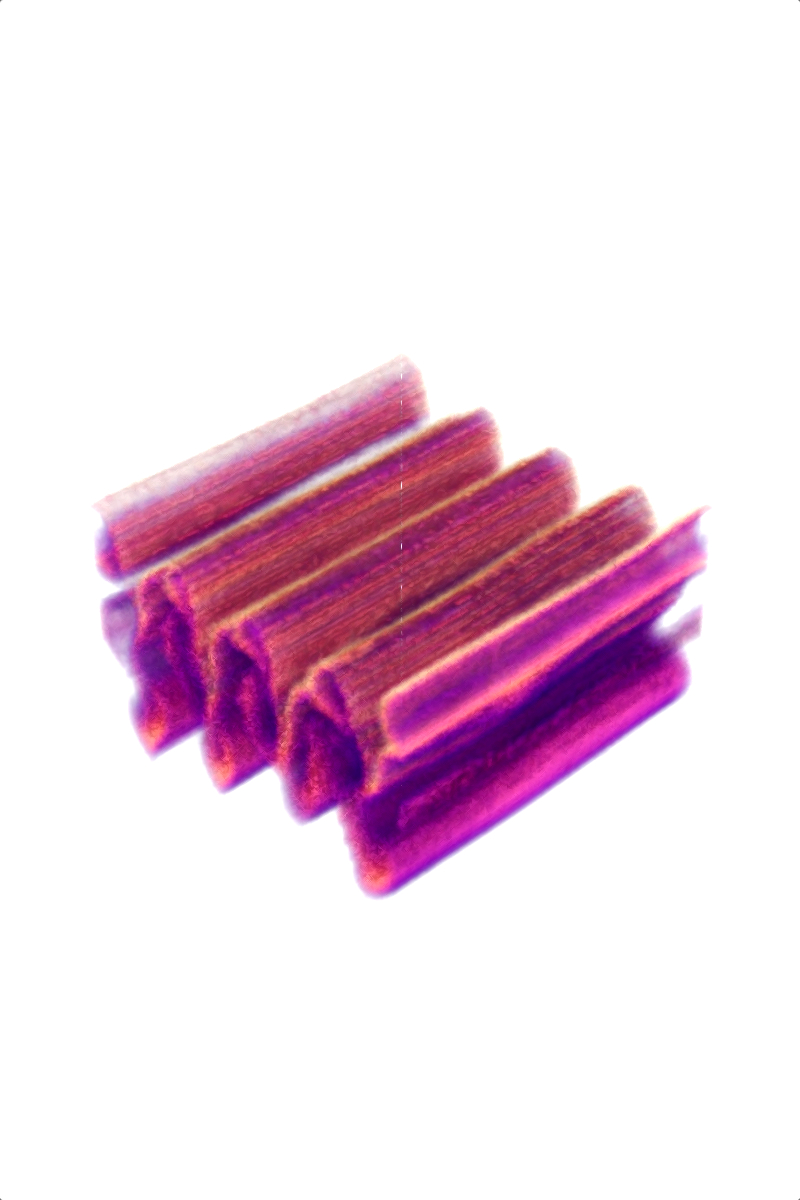}}}
        ~
        \subfloat[$t = 5$]{
        \adjustbox{width=0.24\linewidth,valign=b}{\includegraphics[width=\textwidth]{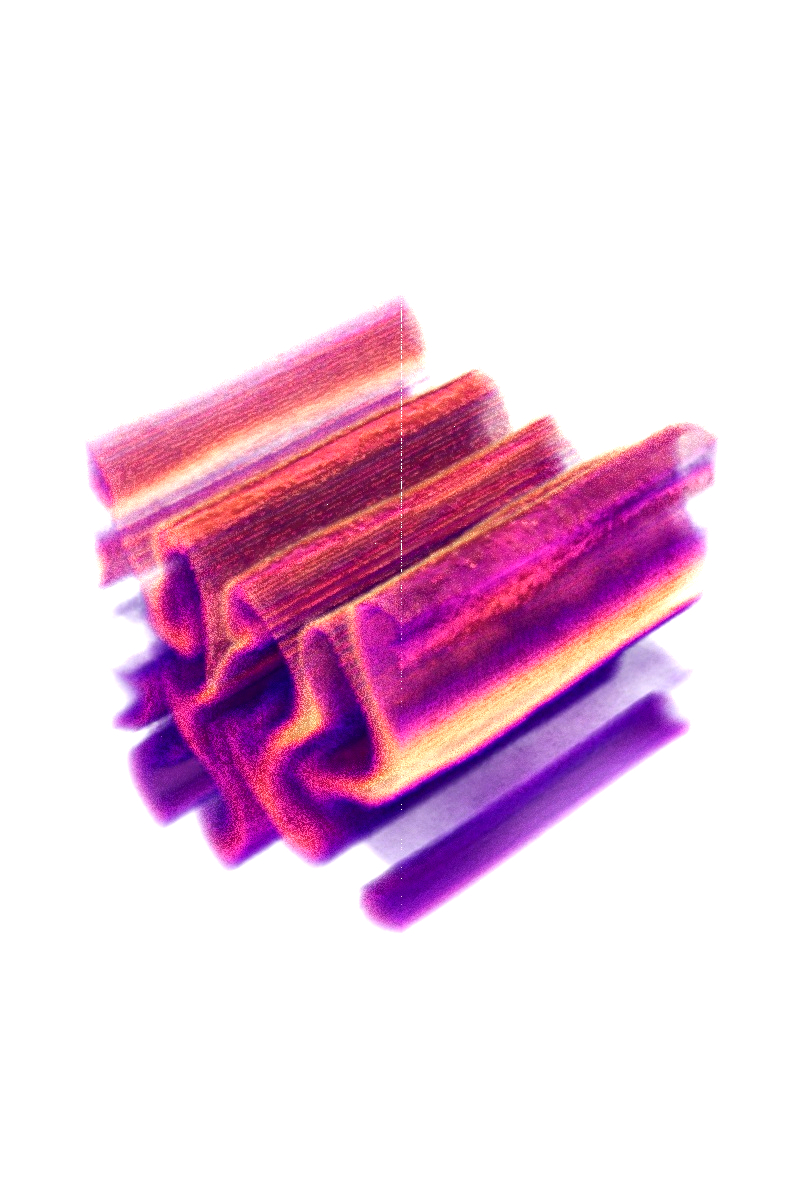}}}
        \newline
        \caption{\label{fig:rt_b7} Volume rendering of the density field for the strongly magnetized ($B_0 = 0.5 B_c$) Rayleigh--Taylor instability problem at varying times computed using a $\mathbb P_3$ scheme on a $N_e = 64 \times 64 \times 128$ mesh.}
    \end{figure}

These three flow conditions were computed using a $\mathbb P_3$ scheme on a $N_e = 64 \times 64 \times 128$ hexahedral mesh with a time step of $\Delta t = 2{\cdot}10^{-5}$. The flow, visualized in the form of a volume rendering of the density field, at various times is shown in \cref{fig:rt_b0}, \cref{fig:rt_b14}, and \cref{fig:rt_b7} for the hydrodynamic, weakly magnetized, and strongly magnetized cases, respectively. For the hydrodynamic case, the canonical flow pattern of the Rayleigh--Taylor instability was observed, with rising bubbles and descending figures. At later times, these features transitioned to a turbulent mixing state. When a weak magnetic field was applied, a significant degree of anisotropy was imparted on the flow, seen in the form of distortions in the bubbles and fingers aligned with the orientation of the magnetic field. Furthermore, the slowed growth rate of the instabilities due to the stabilizing effect of this weak magnetic field could be observed. For the strongly magnetized case, the magnetic field effectively damped all instabilities along the orientation of the field, such that the resulting flow only varied perpendicular to the orientation of the field. Given a long enough simulation time, this flow would be expected to transition to a three-dimensional turbulent mixing state as nonlinear transport effects overcome the stabilizing nature of the magnetic field. For all cases, the flow was numerically well-behaved, indicating that the proposed approach can be effectively applied to three-dimensional flows in both the magnetized and hydrodynamic regimes. 
    
\begin{figure}[tbhp]
    \centering
    \adjustbox{width=0.4\linewidth, valign=b}{\begin{tikzpicture}[spy using outlines={rectangle, height=3cm,width=2.3cm, magnification=3, connect spies}]
\begin{axis} [ybar,
		axis line style={latex-latex},
	    axis x line=left,
        axis y line=left,
        bar width=30pt,
    	xmin=0.5, xmax=2.5,
    	ymin=0, ymax=300,
    	ylabel={GPU hours per characteristic time},
        xtick={1, 2},
        xticklabels={$\mathbb P_3$, $\mathbb P_4$},
        clip mode=individual,
    	legend style={at={(0.03, 0.97)},anchor=north west},
    	legend cell align={left}]
     
\addplot[draw=black,fill=black!80]
    coordinates {
    	(1, 25.3) 
    	(2, 39.8) 
    };
\addlegendentry{Proposed algorithm}

\addplot[draw=black,pattern=north east lines,pattern color = gray]
    coordinates {
    	(0.98, 61.4) 
    	(1.98, 257)
    };
\addlegendentry{Original algorithm}

\node at (0.85, 35.3) {25.3};
\node at (1.15, 71.4) {61.4};
\node at (1.85, 49.8) {39.8};
\node at (2.15, 267) {257};

\end{axis}

\end{tikzpicture}}
    \caption{\label{fig:cost} Comparison of the wall-clock time to reach $t=1$ for the weakly magnetized ($B_0 = 0.1 B_c$) Rayleigh--Taylor instability problem with a $\mathbb P_3$ (left) and $\mathbb P_4$ (right) scheme with the same number of degrees of freedom using the original algorithm of \citet{Dzanic2022} and the proposed algorithm. 
    }
\end{figure}
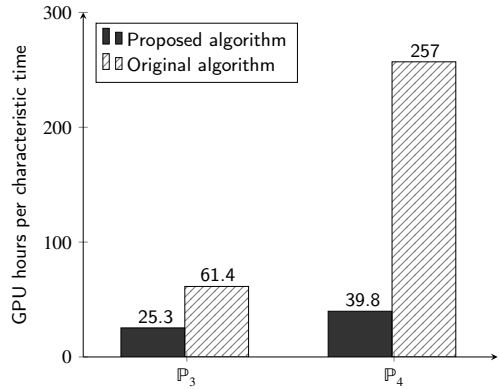

To quantify the efficiency improvements of the proposed algorithm in comparison to the original approach in \citet{Dzanic2022} which utilizes repeated evaluations of the matrix-vector product in \cref{eq:matvec}, a runtime comparison between the two methods was performed for the case of $B_0 = 0.1 B_c$. The cost was evaluated on 16 NVIDIA V100 GPUs with respect to the wall-clock time elapsed until the simulation reached $t=1$, and the results are shown in \cref{fig:cost}. While the original approach required 61.4 GPU hours, the proposed approach required only 25.3 GPU hours, a speedup factor of approximately 2.4 across the entire simulation time. Furthermore, this speedup is expected to increase with higher approximation orders due to the increased number of solution points per element. To confirm this, an identical comparison was performed using a $\mathbb P_4$ approximation on a $N_e = 51 \times 51 \times 102$ mesh, which results in approximately the same number of degrees of freedom. At this approximation order, the original approach required 257 GPU hours whereas the proposed approach required only 39.8 GPU hours, a speedup factor of approximately 6.5. These results indicate that the proposed algorithmic improvements both substantially decrease the overall computational cost of the entropy filtering approach and show much better scaling with respect to the approximation order.

\section{Conclusions}\label{sec:conclusion}
In this work, a positivity-preserving adaptive filtering approach was proposed for shock capturing in discontinuous spectral element approximations of the ideal magnetohydrodynamics equations. The proposed scheme can be considered as an extension of the entropy filtering approach \citep{Dzanic2022} introduced by the authors for the gas dynamics equations to the ideal magnetohydrodynamics system. By formulating convex invariants such as positivity of density and pressure and a local discrete minimum entropy principle as discrete constraints on the solution, the amount of filtering necessary to satisfy the constraints was computed as an element-wise scalar optimization problem. This approach was combined with the eight-wave method of \citet{Powell1999} for enforcing a solenoidal magnetic field. As this method introduced non-conservative source terms to the system, an operator splitting approach was proposed and its effects on the assumptions necessitated by the adaptive filtering approach were analyzed. An improved algorithm for solving the optimization problem for the filter strength was also introduced which significantly improved the computational efficiency of the proposed method. 

The proposed scheme could robustly resolve strong discontinuities while recovering high-order accuracy in smooth regions of the flow and could be easily and efficiently implemented on general unstructured grids. The efficacy of the approach was shown in a variety of numerical experiments, ranging from simple transport and shock tubes to extremely magnetized blast waves and three-dimensional magnetohydrodynamic instabilities. Furthermore, the proposed algorithmic enhancements yielded significant improvements in the computational cost and showed much better scaling with respect to approximation order, reducing the total runtime of the simulations by a factor of 2.4 for $\mathbb P_3$ approximations and 6.5 for $\mathbb P_4$ approximations. Future improvements to the proposed scheme could focus on applying different filter kernels to various components on the solution, alternate methods for enforcing a divergence-free magnetic field, and anisotropic filtering approaches.

\section*{Acknowledgements}
\label{sec:ack}
This work was supported in part by the U.S. Air Force Office of Scientific Research via grant FA9550-21-1-0190 ("Enabling next-generation heterogeneous computing for massively parallel high-order compressible CFD") of the Defense University Research Instrumentation Program (DURIP) under the direction of Dr. Fariba Fahroo.
\bibliographystyle{unsrtnat}
\bibliography{reference}



\end{document}